\documentclass{amsart}

\usepackage{latexsym}
\usepackage{amscd}
\usepackage[all]{xy}
\usepackage[mathscr]{euscript}
\usepackage{dsfont}
\normalsize

\RequirePackage{xspace}

\newcommand{\thought}[1]{}
\renewcommand{\thought}[1]{ \textbf{[#1]}}

\usepackage{enumerate}        
\newenvironment{roenumerate}{\begin{enumerate}[\upshape (i)]}{\end{enumerate}}

\newcommand\nc {\newcommand}
\newcommand\rnc{\renewcommand}

\newcount\blopone
\newcount\xone
\newcount\xtwo
\newcount\ytwo

\newtheorem{theorem}{Theorem}[subsection]
\newtheorem{prop}[theorem]{Proposition}
\newtheorem{axiom}[theorem]{Axiom}
\newtheorem{observation}[theorem]{Observation}
\newtheorem{com}[theorem]{Comment}
\newtheorem{computation}[theorem]{Computation}
\newtheorem{redu}[theorem]{Reduction}
\newtheorem{refinement}[theorem]{Refinement}
\newtheorem{summary}[theorem]{Summary}
\newtheorem{importnota}[theorem]{Important Notation}
\newtheorem{prblm}[theorem]{Problem}
\newtheorem{notation}[theorem]{Notation}
\newtheorem{defin}[theorem]{Definition}
\newtheorem{caution}[theorem]{Caution}
\newtheorem{remark}[theorem]{Remark}
\newtheorem{reminder}[theorem]{Reminder}
\newtheorem{illustration}[theorem]{Illustration}
\newtheorem{strategy}[theorem]{Strategy}
\newtheorem{lemma}[theorem]{Lemma}
\newtheorem{convention}[theorem]{Convention}
\newtheorem{construction}[theorem]{Construction}
\newtheorem{corollary}[theorem]{Corollary}
\newtheorem{example}[theorem]{Example}
\newtheorem{sketch}[theorem]{Sketch}
\newtheorem{conclusion}[theorem]{Conclusion}
\newtheorem{triviality}[theorem]{Triviality}
\newtheorem{proto}[theorem]{Prototype Quasifibration}
\newtheorem{cauex}[theorem]{Cautionary Example}
\newtheorem{hypo}[theorem]{Hypothesis}
\newtheorem{subth}{ }[theorem]
\newtheorem{case}{Case}[theorem]
\newtheorem{ssubth}{ }[subth]
\newtheorem{facts}[theorem]{Facts}
\newtheorem{discussion}[theorem]{Discussion}
\newtheorem{application}[theorem]{Application}

\nc\tri[1]{\begin{triviality}
\label{#1}}
\nc\dis[1]{\begin{discussion}
\label{#1}
\begin{em}}
\nc\apl[1]{\begin{application}
\label{#1}
\begin{em}}
\nc\fac[1]{\begin{facts}
\label{#1}
\begin{em}}
\nc\cas[1]{\begin{case}
\label{#1}
\begin{em}}
\nc\cvn[1]{\begin{convention}
\label{#1}
\begin{em}}
\nc\rfn[1]{\begin{refinement}
\label{#1}}
\nc\prt[1]{\begin{proto}
\label{#1}}
\nc\lem[1]{\begin{lemma}
\label{#1}}
\nc\pro[1]{\begin{prop}
\label{#1}}
\nc\thm[1]{\begin{theorem}
\label{#1}}
\nc\axm[1]{\begin{axiom}
\label{#1}}
\nc\obs[1]{\begin{observation}
\label{#1}}
\nc\cor[1]{\begin{corollary}
\label{#1}}
\nc\dfn[1]{\begin{defin}
\label{#1}}
\nc\sthm[1]{\begin{subth}
\label{#1}}
\nc\exm[1]{\begin{example}
\label{#1}
\begin{em}}
\nc\skt[1]{\begin{sketch}
\label{#1}
\begin{em}}
\nc\plm[1]{\begin{prblm}
\label{#1}
\begin{em}}
\nc\rmk[1]{\begin{remark}
\label{#1}
\begin{em}}
\nc\rmd[1]{\begin{reminder}
\label{#1}
\begin{em}}
\nc\str[1]{\begin{strategy}
\label{#1}
\begin{em}}
\nc\ntn[1]{\begin{notation}
\label{#1}
\begin{em}}
\nc\cmp[1]{\begin{computation}
\label{#1}
\begin{em}}
\nc\smr[1]{\begin{summary}
\label{#1}
\begin{em}}
\nc\cau[1]{\begin{caution}
\label{#1}
\begin{em}}
\nc\hyp[1]{\begin{hypo}
\label{#1}}
\nc\imn[1]{\begin{importnota}
\label{#1}
\begin{em}}
\nc\rdn[1]{\begin{redu}
\label{#1}
\begin{em}}
\nc\cax[1]{\begin{cauex}
\label{#1}
\begin{em}}
\nc\cmt[1]{\begin{com}
\label{#1}
\begin{em}}
\nc\con[1]{\begin{construction}
\label{#1}
\begin{em}}
\nc\ill[1]{\begin{illustration}
\label{#1}
\begin{em}}
\nc\ssthm[1]{\begin{ssubth}
\label{#1}
\begin{em}}
\nc\cnc[1]{\begin{conclusion}
\label{#1}
\begin{em}}

\nc\elem{\end{lemma}}
\nc\erdn{\end{em}\end{redu}}
\nc\erfn{\end{refinement}}
\nc\eprt{\end{proto}}
\nc\ethm{\end{theorem}}
\nc\eaxm{\end{axiom}}
\nc\eobs{\end{observation}}
\nc\ecor{\end{corollary}}
\nc\edfn{\end{defin}}
\nc\esthm{\end{subth}}
\nc\epro{\end{prop}}
\nc\etri{\end{triviality}}
\nc\eexm{\end{em}
\end{example}}
\nc\eskt{\end{em}
\end{sketch}}
\nc\ecvn{\end{em}
\end{convention}}
\nc\ecmt{\end{em}
\end{com}}
\nc\efac{\end{em}
\end{facts}}
\nc\edis{\end{em}
\end{discussion}}
\nc\eapl{\end{em}
\end{application}}
\nc\ermk{\end{em}
\end{remark}}
\nc\ermd{\end{em}
\end{reminder}}
\nc\eill{\end{em}
\end{illustration}}
\nc\eplm{\end{em}
\end{prblm}}
\nc\ecas{\end{em}
\end{case}}
\nc\ecau{\end{em}
\end{caution}}
\nc\ecax{\end{em}
\end{cauex}}
\nc\eimn{\end{em}
\end{importnota}}
\nc\entn{\end{em}
\end{notation}}
\nc\estr{\end{em}
\end{strategy}}
\nc\ecmp{\end{em}
\end{computation}}
\nc\econ{\end{em}
\end{construction}}
\nc\esmr{\end{em}
\end{summary}}
\nc\ehyp{
\end{hypo}}
\nc\ecnc{\end{em}
\end{conclusion}}
\nc\essthm{\end{em}
\end{ssubth}}

\nc\sst{\scriptstyle}
\newcommand{\comment}[1]{}
\newcommand{\ri}{\longrightarrow}
\newcommand{\sr}{\rightarrow}

\newcommand{\zz}{{\mathbb Z}}
\nc\gm{{{\mathbb G}_m}}

\newcommand{\D}{{\mathbf D}}
\newcommand{\Dqc}{{\mathbf D_{\text{\bf qc}}}}

\newcommand{\pp}{{\mathbb P}}

\newcommand{\oo}{\otimes}

\newcommand{\ak}{{\mathbb A}}
\newcommand\cpl{\mathbb{C}}
\nc\op{^{\hbox{\rm\tiny op}}}
\nc\mth{^{\hbox{\rm\tiny th}}}

\nc\dcoh{\D_{\mathbf{coh}}^b}
\nc\dcohm{\D_{\mathbf{coh}}^-}
\nc\dcohp{\D_{\mathbf{coh}}^+}
\nc\script{\mathscr}
\nc\z{\zeta}
\nc\bc{{\mathbb{BC}}}
\nc\ct{{\script T}}
\nc\cf{{\script F}}
\nc\cl{{\script L}}
\nc\cv{{\script V}}
\nc\cq{{\script Q}}
\nc\cu{{\script U}}
\nc\ce{{\script E}}
\nc\cg{{\script G}}
\nc\ch{{\script H}}
\nc\cs{{\script S}}
\nc\car{{\script R}}
\nc\cd{{\script D}}
\nc\cc{{\script C}}
\nc\ck{{\script K}}
\nc\ca{{\script A}}
\nc\ci{{\script I}}
\nc\cj{{\script J}}
\nc\co{{\script O}}
\nc\cm{{\script M}}
\nc\cz{{\script Z}}
\nc\cx{{\script X}}
\nc\cy{{\script Y}}
\nc\cw{{\script W}}
\nc\bd{\begin{description}}
\nc\ed{\end{description}}
\nc\ctob{{\script C}at\big(\ci^{op},\ca\big)}
\nc\clim{{\ds\mathop{\rm lim}_{\ds\longleftarrow}}\,}
\nc\climi{\clim^{\!i}\,}
\nc\climn{\clim^{\!n}\,}
\nc\colim{{\ds\mathop{\rm colim}_{\ds\la}}}
\nc\oa{\overline{\ca}}
\nc\s{\sigma}
\nc\ta{\tau}
\nc\os{\overline\sigma}
\nc\ot{\overline\tau}
\nc\T{\Sigma}
\nc\Tm{\Sigma^{-1}}
\nc\de[1]{{\mathop{\rm deg(#1)}}}
\nc\Ad[1]{\mathop{\rm Ad}(#1)}
\nc\ad[1]{\mathop{\rm ad}(#1)}
\nc\kth{{\it K}--theory}
\nc\loc[1]{{\text{\rm Loc(#1)}}}
\nc\coloc[1]{{\text{\rm Coloc}(#1)}}

\nc\one{\mathds{1}}

\def\der #1 {D\left(#1\right)}
\nc\prf{\begin{proof}}
\nc\eprf{\end{proof}}
\nc\ds{\displaystyle}
\nc\Tor{\text{\rm Tor}}

\nc\cb{{\script B}}
\nc\ab{{\script A}b}

\nc\be{\begin{roenumerate}}
\nc\ee{\end{roenumerate}}

\nc\cat[1]{{\script C}at\Big({\big\{#1\big\}}\op\,\,,\,\,\ab\Big)}
\nc\csab{{\script C}at\big(\cs^{op},\ab\big)}
\nc\ctab{{\script C}at\Big({\{\ct^\alpha\}}^{op},\ab\Big)}
\nc\csex{{\script E}x\big(\cs^{op},\ab\big)}
\nc\ctex{{\script E}x\Big({\{\ct^\alpha\}}^{op},\ab\Big)}
\nc\sub{\qquad\subset\qquad}
\nc\ctr[1]{{\left.\ct\left(-,#1\right)\right|}_{\cs}}
\nc\ctrf[2]{{\left.\ct\left(#1,#2\right)\right|}_{\cs}}
\nc\Ctr[1]{{\left.\ct\left(-,#1\right)\right|}_{\ct^\alpha}}
\nc\Ctrf[2]{{\left.\ct\left(#1,#2\right)\right|}_{\ct^\alpha}}

\nc\la{\longrightarrow}
\nc\nin{\noindent}
\nc\cad[1]{\text{card}(#1)}
\nc\eq{\quad=\quad}
\nc\BA{\begin{array}{c}}
\nc\EA{\end{array}}
\nc\barr{
\[
\begin{array}{cccccccccccccccc}
}
\nc\earr{
\end{array}
\]
}
\nc\as[1]{{\langle S\rangle}^{#1}}
\nc\sh{\text{\it shift}}

\nc\yy[1]{{\left.\ct\left(-,#1\right)\right|}_{\ct^c}}
\nc\vrep[2]{{\left.\ct\left(#1,#2\right)\right|}_{\ct^\alpha}}
\nc\da{\downarrow}
\nc\Hom{{\mathop{\rm Hom}}}
\nc\RHom{{\mathbf{R}\text{\rm Hom}}}
\nc\HHom{{\script H}{\mathop{\rm om}}}
\nc\RHHom{{\script{RH}}{\mathop{\rm om}}}
\nc\End{{\mathop{\rm End}}}
\nc\Ext{{\mathop{\rm Ext}}}
\nc\mr\Modtc
\nc\PExt{{\mathop{\rm PExt}}}
\nc\stm{\text{\rm stmod}(kG)}
\nc\stM{\text{\rm StMod}(kG)}
\nc\e{\varepsilon}
\nc\p{\mathfrak{p}}

\nc\rs{\s^{-1}A}
\nc\br{{\{\s^{-1}A\}}}
\nc\ra\ri
\nc\y[1]{\mathbf{y}#1}
\nc\x[1]{\mathbf{z}#1}
\nc\mmod[1]{#1\text{--\rm mod}}
\nc\Mod[1]{#1\text{--\rm Mod}}
\nc\Md {\ensuremath{\mathop{\textup{Mod}}}}
\rnc\mod[1]{\ensuremath{\mathop{\textup{mod-}#1}}\xspace}
\nc\Modtc{\Mod{\ct^c}}
\nc\pgldim[1]{\mathop{\rm pgldim}\,#1}
\nc\tf{{\rm [TR5]}}
\nc\tfs{{\rm [TR5$^*$]}}
\nc\Fun{\text{\rm Funct}(F\op,\ab)}
\nc\sym{\text{\rm Sym}}
\nc\sgn{\text{\rm sgn}}
\nc\Pro{\text{\rm Prod}^{}_\alpha(F\op,\ab)}
\nc\Yt[1]{{\left.\Hom_\ct^{}\left(-,#1\right)\right|}_F^{}}
\nc\dl{\delta}
\nc\Proj[1]{#1\text{--\rm Proj}}
\nc\proj[1]{#1\text{--\rm proj}}
\nc\Flat[1]{\text{\rm Flat}\,#1}
\nc\Inj[1]{\text{\rm Inj\,}#1}
\nc\qc[1]{\text{\rm qc\,}#1}
\nc\ov{\overline}
\nc\wt{\widetilde}
\nc\wh{\widehat}
\nc\ph{\varphi}
\nc\tstr{{\it t}--structure}
\nc\spec[1]{{\text{\rm Spec}(#1)}}
\nc\Spec[1]{{\text{\rm Spec}\big(#1\big)}}
\newcommand{\m}{\mathfrak{m}}
\newcommand{\fC}{\mathfrak{C}}

\newcommand{\Dqcpl}{{\mathbf D^+_{\mathbf{qc}}}}

\nc\EProd{\text{\rm EProd}}
\nc\ECoprod{\text{\rm ECoprod}}
\nc\Prod{\text{\rm Prod}}
\nc\Coprod{\text{\rm Coprod}}
\nc\COprod{\text{\rm coprod}}
\nc\ldimp{\text{\rm LDim}^{\prod}}
\nc\ldimc{\text{\rm LDim}^{\coprod}}

\nc\hoco{
\begin{picture}(40,10)
\put(20,0){\makebox(0,0)[b]{\text{\rm Hocolim}}}
\put(5,-2){\vector(1,0){30}}
\end{picture}\,}

\nc\holim{
\begin{picture}(40,10)
\put(20,0){\makebox(0,0)[b]{\text{\rm Holim}}}
\put(35,-2){\vector(-1,0){30}}
\end{picture}}

\nc\Cop{\text{\rm Coprod}}
\nc\seq{{\mathbb{S}_{\mathbf{e}}}}
\nc\se{{S^{\text{\tt e}}}}
\nc\te{{T^{\text{\tt e}}}}
\nc\LL{{\text{\bf L}}}
\nc\R{\text{\bf R}}
\nc\id{\text{\rm id}}
\nc\ev{\text{\rm ev}}
\nc\supp{\text{\rm supp}}
\nc\Loc{\text{\rm Loc}}
\nc\Thick{\text{\rm Thick}}

\begin{document}

\author{Amnon Neeman}\thanks{The research was partly supported 
by the Australian Research Council}
\address{Centre for Mathematics and its Applications \\
        Mathematical Sciences Institute\\
        Building 145\\
        The Australian National University\\
        Canberra, ACT 2601\\
        AUSTRALIA}
\email{Amnon.Neeman@anu.edu.au}

\title{Grothendieck duality made simple}

\begin{abstract}
  It has long been accepted that the foundations of Grothendieck duality
  are complicated. This has changed recently.

  By ``Grothendieck duality'' we mean what, in the old literature, used to go
  by the name ``coherent duality''. This isn't to be confused with
  what is nowadays called ``Verdier duality'', and
  used to pass as ``$\ell$-adic duality''. 
  The footnote below comments
  on the historical inaccuracy of the modern terminology.
\end{abstract}\footnote{
  The prevailing current
  terminology---for duality in \'etale cohomology,
  that is ``$\ell$-adic duality''---is
  historically incorrect. The idea was originally due
  not to Verdier but to Grothendieck, see
  his work in SGA5 on what is nowadays called the
  formalism of the six operations. Since this survey is
  about coherent duality we elaborate no further.}



\maketitle

\tableofcontents

\section{Introduction}
\label{S-1}

There are two classical paths to the foundations of Grothendieck
duality: one due to Grothendieck and Hartshorne~\cite{Hartshorne66} and
(much later)
Conrad~\cite{Conrad00}, and a second due to
Deligne~\cite{Deligne66} and 
Verdier~\cite{Verdier68} and (much later)
Lipman~\cite{Lipman09}. The consensus has been that both are
unsatisfactory. If you listen to the detractors of the respective
approaches: the first is a nightmare to set up, the second leads to
a theory where you can't compute anything. While exaggerated, each criticism
used to have a kernel of truth to it. Lipman summed
it up more circumspectly and diplomatically some years ago,
by saying there is no royal road to the subject.

And Lipman is probably the person who worked hardest on simplifying
the foundations.

In passing let us mention that, while the
two foundational avenues to setting up the subject must
obviously be related, the details of this link are far
from clear---in fact they haven't yet
been fully worked out.
And Lipman also happens to be the person who has tried hardest to understand
this bridge.

Back to the history of
the field: what happened is that in the mid 1970s the math
community gradually started losing interest in the
project, and by the mid 1980s the exodus was
all but complete---most people had given up on improving the
foundations and moved
on. For three decades now there have been at most a dozen people actively
working in the field---they break up into two groups: Lipman, his
students and collaborators and Yekutieli, his students and collaborators.
These two groups made regular, incremental progress, see for example
Alonso, Jerem\'{i}as and Lipman~\cite{Alonso-Jeremias-Lipman97,Alonso-Jeremias-Lipman99,Alonso-Jeremias-Lipman11,Alonso-Jeremias-Lipman12},
Lipman~\cite{Lipman84,Lipman87,Lipman02},
Lipman, Nayak and Sastry~\cite{Lipman-Nayak-Sastry05},
Lipman and Neeman~\cite{Lipman-Neeman07},
Nayak~\cite{Nayak05,Nayak09},
Porta, Shaul and Yekutieli~\cite{Porta-Shaul-Yekutieli14},
Shaul~\cite{Shaul15,Shaul16,Shaul17},
Sastry~\cite{Sastry04,Sastry05},
Yekutieli~\cite{Yekutieli92,Yekutieli92A,Yekutieli10,Yekutieli16,Yekutieli19} and 
Yekutieli and Zhang~\cite{Yekutieli-Zhang97,Yekutieli-Zhang99,Yekutieli-Zhang05,
  Yekutieli-Zhang06,Yekutieli-Zhang08,Yekutieli-Zhang09}.

And then, a couple of years ago, there was a seismic shift. In this article
we attempt to describe the recent progress in a way accessible to the
non-expert.

Let us stress that the term ``non-expert'' in the last paragraph
is to be understood in the strong sense: the intended audience of
this survey includes non-algebraic-geometers, 
and familiarity with derived categories isn't assumed.
The reason for the inclusive
exposition is that the
questions opened up by the recent progress might well interest
people in diverse fields, the most obvious being Hochschild homology and
cohomology. The Hochschild experts might wish to start with
Computation~\ref{C30.41}, as well as Problems~\ref{P95.7}, \ref{P95.999} and
\ref{P95.3.119}. Computation~\ref{C30.41} spells out exactly where
and how Hochschild 
homology and cohomology played a key role in the breakthrough
we report.
The three problems suggest obvious variants, explaining why each would be 
interesting to solve. To put it in a nutshell:
up to the present time we---meaning
the handful of algebraic geometers still working on the foundations
of Grothendieck duality---have only been
able to carry out the Hochschild homology computation of~\ref{C30.41}, which
amounts to a very simple, special case of the general problem.
Given how profoundly this baby computation has transformed our
understanding of Grothendieck duality, we
warmly invite mathematicians more adept and dextrous at using
the Hochschild machinery to come to our aid.

\medskip

\nin{\bf Acknowledgements.}\ \
The author would like to thank Leo Alonso,
Asilata Bapat, Spencer Bloch,
Jim Borger, Jesse Burke,
Pierre Deligne, Anand Deopurkar,
Luc Illusie, Ana Jerem\'{\i}as,
Steve Lack, Joe Lipman, Bregje Pauwels, Pramath Sastry, Ross Street,
Michel Van den Bergh, Amnon Yekutieli and an anonymous referee
for questions, comments,
corrections and improvements---based both on earlier versions
of this manuscript, and on talks I've given presenting
parts or all of
the material. The contributions of the people listed have
greatly influenced the
exposition. Needless to say the flaws and mistakes that remain are
entirely my fault.

\section{Background}
\label{S24}

\subsection{Reminder: formally inverting morphisms}
\label{SS24.257}

Let $\cc$ be a category, and let $S\subset\text{Mor}(\cc)$ be some
collection of morphisms. It is a theorem of
Gabriel and Zisman~\cite{Gabriel-Zisman67} that one may form a functor
$F:\cc\la S^{-1}\cc$ so that
\be
\item
  The functor $F$ takes every element of $S$ to an isomorphism.
\item
If $H:\cc\la\cb$ is a functor, taking every element of $S$
to an isomorphism, then there exists a unique functor $G:S^{-1}\cc\la\cb$
rendering commutative the triangle
\[\xymatrix@C+40pt@R-20pt{
        & S^{-1}\cc\ar@{.>}[dd]^-{\exists!G} \\
\cc\ar[ru]^-F\ar[rd]_-H &  \\
              & \cb
}\]
\ee
We call this construction \emph{formally inverting the morphisms in $S$.}

\rmk{R24.2579}
On objects the functor $F$ is the identity: the objects of
$S^{-1}\cc$ are identical to those of $\cc$. But the morphisms in $S^{-1}\cc$
are complicated. Clearly any morphism of $\cc$ must have an image
in $S^{-1}\cc$, but $S^{-1}\cc$ must also contain inverses of the
images of morphisms in $S$. And then we must be able to compose any
finite string of these.

The morphisms of $S^{-1}\cc$ are in fact equivalence classes of such
finite strings. The problem becomes to figure out when two such
strings are equivalent, that is which strings must have the same
composite in $S^{-1}\cc$. This is usually called
\emph{the calculus of fractions} of $S^{-1}\cc$. And without an
understanding of this calculus of fractions the category $S^{-1}\cc$
is unwieldy.

The category $S^{-1}\cc$ can be dreadful in general, for example:
it may happen that $\cc$ has small Hom-sets but $S^{-1}\cc$ doesn't.
\ermk

\subsection{Reminder: derived categories and derived functors}
\label{SS24.250}
Let $\ca$
be an abelian category. We will be looking at categories we will denote
$\D_{\fC}^{}(\ca)$. The category $\D_{\fC}^{}(\ca)$ is as follows:
\be
\item
The objects: an object in $\D_\fC^{}(\ca)$ is a cochain complex of
objects in $\ca$, that is a diagram in $\ca$
\[\xymatrix@C+5pt{
\cdots \ar[r]& A^{-2}  \ar[r]& A^{-1}  \ar[r]& A^{0}  \ar[r]
&  A^{1} \ar[r]& A^{2} \ar[r] &\cdots
}\]
where the composite $A^i\la A^{i+1}\la A^{i+2}$ vanishes for
every $i\in\zz$. The $\fC$ in
$\D_\fC^{}(\ca)$ stands for conditions: we may impose conditions on the objects,
it may be that not every cochain complex belongs to $\D_\fC^{}(\ca)$.
\item
Morphisms: cochain maps are examples, that is commutative diagrams
\[\xymatrix@C+5pt{
\cdots \ar[r]& A^{-2}  \ar[d]\ar[r]& A^{-1} \ar[d] \ar[r]& A^{0}\ar[d]  \ar[r]
&  A^{1} \ar[d]\ar[r]& A^{2}\ar[d] \ar[r] &\cdots\\
\cdots \ar[r]& B^{-2}  \ar[r]& B^{-1}  \ar[r]& B^{0}  \ar[r]
&  B^{1} \ar[r]& B^{2} \ar[r] &\cdots
}\]
where the rows are objects in $\D_\fC^{}(\ca)$.
But we also formally invert the cohomology isomorphisms.
\ee

\rmk{R24.989796}
In the special case of $\D_\fC^{}(\ca)$
the calculus of fractions is reasonably well understood,
there's a rich literature about it---but we will not explain this here. This
means that, whenever we tell the reader about some computation
of morphisms, we will be asking the beginner to accept it on faith.
The expert will notice that all the computations we mention are
easy.
\ermk

\exm{E24.1}
Let $R$ be a commutative ring,
and let $\ca=\Mod R$ be the abelian category of $R$--modules.
Then the category $\D(\Mod R)$ has for its objects all the cochain
complexes of $R$--modules. The category $\D_{K\text{--Flat}}^{}(\Mod R)$
has for its objects all the $K$-flat complexes, the category
$\D_{K\text{--Inj}}^{}(\Mod R)$ has for its objects all the $K$--injective
complexes.

We remind the reader: a complex $F^*$ is $K$--flat if, for all
acyclic complexes $A^*$, the complex $A^*\oo_R^{} F^*$ is acyclic. A complex
$I^*$ is $K$--injective if, for all
acyclic complexes $A^*$, the complex $\Hom_R^{}(A^*,I^*)$ is
acyclic.
\eexm

Next we recall functors. If $F:\ca\la\cb$ is an additive functor, we often
want to extend $F$ to derived categories. But the simple-minded
approach does not in general work, you cannot simply apply $F$ to the
cochain complexes.

\exm{E24.3}
Let $f:R\la S$ be a homomorphism of commutative rings,
let $\ca=\Mod R$ and let $\cb=\Mod S$.
Fix an object $A\in\Mod S$.
We wish to consider the functor $A\oo_R^{}(-)$, that is the functor that
takes the object $B\in\Mod R$ to the object $A\oo_R^{} B$ in $\Mod S$.

If we try to extend it to a functor
$A\oo_R^{}(-):\D(\Mod R)\la\D(\Mod S)$ we run
into the following problem: it is entirely possible to have a
cochain map $\wt B^*\la B^*$ of $R$--modules, inducing an isomorphism in
cohomology, but where $A\oo_R^{} \wt B^*\la A\oo_R^{} B^*$
is \emph{not} a cohomology
isomorphism. The special case where $\wt B^*=0$ is already problematic.
A cochain map
\[\xymatrix@C+5pt{
\cdots \ar[r]& 0  \ar[d]\ar[r]& 0 \ar[d] \ar[r]& 0\ar[d]  \ar[r]
&  0 \ar[d]\ar[r]& 0\ar[d] \ar[r] &\cdots\\
\cdots \ar[r]& B^{-2}  \ar[r]& B^{-1}  \ar[r]& B^{0}  \ar[r]
&  B^{1} \ar[r]& B^{2} \ar[r] &\cdots
}\]
is a cohomology isomorphism if $B^*$ is acyclic. But without some restrictions
we would not expect $A\oo_R^{} B^*$ to be acyclic, meaning
\[\xymatrix@C-3pt{
\cdots \ar[r]& 0  \ar[d]\ar[r]& 0 \ar[d] \ar[r]& 0\ar[d]  \ar[r]
&  0 \ar[d]\ar[r]& 0\ar[d] \ar[r] &\cdots\\
\cdots \ar[r]& A\oo_R^{} B^{-2}  \ar[r]& A\oo_R^{} B^{-1}  \ar[r]& A\oo_R^{} B^{0}  \ar[r]
&  A\oo_R^{} B^{1} \ar[r]& A\oo_R^{} B^{2} \ar[r] &\cdots
}\]
will not be a cohomology isomorphism.
\eexm

\con{C24.5}
As above, suppose we are given an additive functor $F:\ca\la\cb$.
The remedy is to pass to ``derived functors''. The idea is as follows:
\be
\item
Find a condition $\fC$ on the objects of $\D(\ca)$, so that
if a cochain map $A^*\la B^*$ between objects in $\D_\fC(\ca)$
is a cohomology isomorphism then so is
$FA^*\la FB^*$.
\item
Prove that the natural functor $I:\D_\fC(\ca)\la\D(\ca)$ is an equivalence of
categories.
\ee
Once we achieve (i) and (ii) above, we declare the derived functor of
$F$ to be the composite
\[\xymatrix@C+15pt{
\D(\ca) \ar[r]^-{I^{-1}} &\D_\fC(\ca)\ar[rr]^-F &&\D(B)
}\]
\econ

\exm{E24.7}
Let us return to the situation of Example~\ref{E24.3}: we are given a
ring homomorphism $R\la S$ and an $S$--module $A$.
Then, while the functor $F(-)=A\oo_R^{}(-)$ does not respect general cochain
maps inducing cohomology isomorphisms, it does respect them when
the cochain complexes are $K$--flat as in Example~\ref{E24.1}.
It turns out that the natural functor
$I:\D_{K\text{--Flat}}^{}(\Mod R)\la\D(\Mod R)$
is an equivalence of categories, and we define the functor
$A\oo^\LL_{R}(-):\D(\Mod R)\la\D(\Mod S)$ to be the composite
\[\xymatrix@C+10pt{
\D(\Mod R) \ar[r]^-{I^{-1}} &\D_{K\text{--Flat}}^{}(\Mod R)\ar[rr]^-{A\oo_R^{}(-)} &&\D(\Mod S)
}\]

Now consider the functor
$\Hom_R^{}(A,-):\Mod R\la\Mod S$. Once again this functor does not respect
general cochain maps inducing cohomology isomorphisms. But it does
respect them if the cochain complexes are $K$--injective as in
Example~\ref{E24.1}, and the natural functor
$I:\D_{K\text{--Inj}}^{}(\Mod R)\la\D(\Mod R)$
is an equivalence of categories. We define the functor
$\RHom_R^{}(A,-):\D(\Mod R)\la\D(\Mod S)$ to be the composite
\[\xymatrix@C+15pt{
\D(\Mod R) \ar[r]^-{I^{-1}} &\D_{K\text{--Inj}}^{}(\Mod R)\ar[rr]^-{\Hom_R^{}(A,-)} &&\D(\Mod S)
}\]
\eexm

\rmk{R24.500}
If $\ca$ is an abelian category, the category $\cc(\ca)$ has the same objects
as $\D(\ca)$ but the only morphisms of $\cc(\ca)$ are the genuine cochain maps.
Inverses of cohomology isomorphisms are not allowed.

Generalizing the discussion of
Example~\ref{E24.7}, we will allow ourselves to derive additive
functors $\cc(\ca)\la\cc(\cb)$. For example: if $A^*$ is an object of
$\cc(\Mod S)$ there are standard functors
\[
A^*\oo_R^{}(-):\cc(\Mod R)\la\cc(\Mod S)\ ,\qquad
\Hom_R^{}(A^*,-):\cc(\Mod R)\la\cc(\Mod S)\ .
\]
The derived functors $A^*\oo^\LL_{R}(-)$ and $\RHom_R^{}(A^*,-)$ are,
respectively, the
composites
\[\xymatrix@R-20pt@C+10pt{
\D(\Mod R) \ar[r]^-{I^{-1}} &\D_{K\text{--Flat}}^{}(\Mod R)\ar[rr]^-{A^*\oo_R^{}(-)} &&\D(\Mod S)\\
\D(\Mod R) \ar[r]^-{I^{-1}} &\D_{K\text{--Inj}}^{}(\Mod R)\ar[rr]^-{\Hom_R^{}(A^*,-)} &&\D(\Mod S)
}\]
\ermk

\rmk{R24.567489}
As  presented in Example~\ref{E24.7} and Remark~\ref{R24.500}
the construction involves
an arbitrary choice. More precisely: in Remark~\ref{R24.500} the 
functor $A^*\oo^\LL_R(-)$ was defined by observing
\be
\item
The natural functor
$I:\D_{K\text{--Flat}}^{}(\Mod R)\la\D(\Mod R)$
is an equivalence of categories.
\item
On the category $\D_{K\text{--Flat}}^{}(\Mod R)$ the functor
$A^*\oo_R^{}(-)$ is well-defined in the obvious way, meaning that when
we restrict to $K$--flat cochain complexes
the classical functor $A^*\oo_R^{}(-)$
respects cochain maps inducing cohomology isomorphisms.
\setcounter{enumiv}{\value{enumi}}
\ee
This allowed us to form  the functor $A^*\oo^\LL_R(-)$ as the composite
\[\xymatrix@C+10pt{
\D(\Mod R) \ar[r]^-{I^{-1}} &\D_{K\text{--Flat}}^{}(\Mod R)\ar[rr]^-{A^*\oo_R^{}(-)} &&\D(\Mod S)
}\]
But the observant reader will note that 
\be
\setcounter{enumi}{\value{enumiv}}
\item
The natural functor
$I:\D_{K\text{--Inj}}^{}(\Mod R)\la\D(\Mod R)$
is also an equivalence of categories.
\item
On the category $\D_{K\text{--Inj}}^{}(\Mod R)$ the functor
$A^*\oo_R^{}(-)$ is also well-defined in the obvious way.
\ee
Hence there is nothing to stop us from forming the composite
\[\xymatrix@R-20pt@C+10pt{
\D(\Mod R) \ar[r]^-{I^{-1}} &\D_{K\text{--Inj}}^{}(\Mod R)\ar[rr]^-{A^*\oo_R^{}(-)} &&\D(\Mod S)
}\]
And it turns out that the composite $A^*\oo^\LL_R(-)$, defined
using (i) and (ii), \emph{does not in general agree} with the
composite defined using (iii) and (iv).

For a choice-free description one notes that, with the category
$\cc(\Mod R)$ as in Remark~\ref{R24.500} and with
$F:\cc(\Mod R)\la D(\Mod R)$ the Gabriel-Zisman quotient map
of Reminder~\ref{SS24.257}, we have a triangle
\[\xymatrix@C-10pt@R-20pt{
 & & & & &\D(\Mod R)\ar[dddd]^-{A^*\oo_R^\LL(-)}\\
 && &\, \ar@/^0.5pc/@{=>}[dd]^\gamma &  & \\
\cc(\Mod R)\ar[uurrrrr]^-{F}\ar[ddrrrrr]_-{A^*\oo_R^{}(-)} & & & & & \\
 & & &\, & & \\
 & & & & & \D(\Mod S)
}\] 
That is: there is a natural transformation from the
composite functor $[A^*\oo_R^\LL(-)]\circ F$ to the functor
$A^*\oo_R^{}(-)$.
And what turns out to be true is that the triangle above
has the obvious universal property: for any triangle 
\[\xymatrix@C-10pt@R-20pt{
 & & & & &\D(\Mod R)\ar[dddd]^-{G}\\
 && &\, \ar@/^0.5pc/@{=>}[dd]^\delta &  & \\
\cc(\Mod R)\ar[uurrrrr]^-{F}\ar[ddrrrrr]_-{A^*\oo_R^{}(-)} & & & & & \\
 & & &\, & & \\
 & & & & & \D(\Mod S)
}\]
of functors and natural tranformation, there is a unique
natural transformation $\ph:G\Longrightarrow[A^*\oo_R^\LL(-)]$
such that $\delta=\gamma\circ(\ph F)$. In category theoretic
language: \emph{the functor $A^*\oo_R^\LL(-)$ is the right Kan
extension of $A^*\oo_R^{}(-)$ along $F$.} 

This description as the right Kan extension is what earns the functor
$A^*\oo_R^\LL(-)$ the title of the \emph{left derived functor} of
the functor $A^*\oo_R^{}(-)$, and is the reason for the $\LL$
in the symbol.
The functor $\RHom(A^*,-)$, which we met 
in  Example~\ref{E24.7} and Remark~\ref{R24.500},
turns out to be the \emph{left Kan extension,} meaning the triangle
of functors and natural transformation is
\[\xymatrix@C-10pt@R-20pt{
 & & & & &\D(\Mod R)\ar[dddd]^-{\RHom(A^*,-)}\\
 && &\,  &  & \\
\cc(\Mod R)\ar[uurrrrr]^-{F}\ar[ddrrrrr]_-{\Hom(A^*,-)} & & & & & \\
 & & &\,\ar@/_0.5pc/@{=>}[uu]_\gamma & & \\
 & & & & & \D(\Mod S)
}\] 
and the universality is with respect to all such triangles. Being the
left Kan extension 
earns the functor $\RHom(A^*,-)$ the title of \emph{right
derived functor of $\Hom(A^*,-)$,} as well as the $\R$ in the
symbol.
\ermk

\subsection{Conventions}
\label{SS24.379}

Unless otherwise stated
all rings are assumed commutative and noetherian, all schemes are
assumed noetherian, and all morphisms of schemes are assumed of finite
type.
Since we will often deal with the derived category $\D(\Mod R)$
we abbreviate it to $\D(R)$. If $X$ is a
scheme we will use $\Dqc(X)$ as a shorthand
for the category $\D_{\mathbf{qc}}^{}(\Mod{\co_X^{}})$. That
is: the objects are cochain complexes of sheaves of $\co_X^{}$--modules,
and the condition we impose is that the cohomology
sheaves are quasicoherent.

\section{Statements of the main results}
\label{S29}

\subsection{Generalities}
\label{SS-1.1}
We begin by setting up the framework.

\rmd{R-1.-1}
Suppose $f:X\la Y$ is a morphism of schemes. There are three induced
functors\footnote{Algebraic geometers might find the symbol $f^\times$
  unfamiliar; the pre-2009 literature on Grothendieck duality
  talks almost exclusively about another functor $f^!$. The functors
  $f^\times$ and $f^!$ agree when
  $f$ is proper, but not in general.
  There is a discussion of $f^!$ and its relation with
  $f^\times$ in Reminder~\ref{R95.-77}, and a brief summary of the
  history in Remark~\ref{R95.-1}. Until we reach that point, in this paper 
we will work exclusively with $f^\times$.}
on the derived categories
\[\xymatrix@R+15pt{
\Dqc(X)\ar[d]|-{\R f_*} \\
\Dqc(Y)\ar@<3ex>[u]^-{\LL f^*} \ar@<-3ex>[u]_-{f^\times}
}\]
where each functor is left adjoint to the one to its right;
in category theoretic notation we write $\LL f^*\dashv\R f_*\dashv f^\times$.
We remind the 
reader what these functors do.
\be
\item
  The functor $\LL f^*$ is the left-derived pullback functor.
  We compute it as in Construction~\ref{C24.5}: let
  $\D_{\mathrm{qc},K\text{\rm--Flat}}^{}(Y)$ be the derived
  category of
  complexes of $\co_Y^{}$--modules, which are $K$--flat
  and have quasicoherent cohomology, and
  let $I:\D_{\mathrm{qc},K\text{\rm--Flat}}^{}(Y)\la\Dqc(Y)$ be the natural map.
  The functor $I$ happens to be an equivalence. 
To evaluate $\LL f^*$ on an object 
$C\in\Dqc(Y)$ you first form $I^{-1}(C)\in\D_{\mathrm{qc},K\text{\rm--Flat}}^{}(Y)$,
then pull back to
obtain the complex $f^{-1}I^{-1}(C)$
  on $X$, and finally form on $X$ the tensor
  product $\LL f^*C=\co_X^{}\oo_{f^{-1}\co_Y^{}}^{}f^{-1}I^{-1}(C)$.
\item
  The functor $\R f_*$ is the right-derived pushforward functor.
  Once again we compute it as in Construction~\ref{C24.5}: this time let
  $\D_{\mathrm{qc},K\text{\rm--Inj}}^{}(X)$ be the derived category of
  complexes of $\co_X^{}$--modules, which are $K$--injective
  and have quasicoherent cohomology, and
  let $I:\D_{\mathrm{qc},K\text{\rm--Inj}}^{}(X)\la\Dqc(X)$ be the natural map.
  This functor $I$ also happens to be an equivalence. 
To evaluate
$\R f_*$ on an object $D\in\Dqc(X)$ you first form
$I^{-1}(D)\in\D_{\mathrm{qc},K\text{\rm--Inj}}^{}(X)$, and then push forward to obtain the complex $\R f_*D=f_*I^{-1}(D)$
  on $Y$.
\item
  The functor $f^\times$ is the mysterious one, Grothendieck
  duality is about understanding its properties. For a general $f$ it 
turns out that $f^\times$ need not be the derived functor of any
functor on abelian categories, with the 
notation as in Remark~\ref{R24.567489}.
The reader might wish to look at Remark~\ref{R95.3.1} and at 
Appendix~\ref{A97} for further discussion of how 
unusual $f^\times$ can be.
\ee
\ermd

\dis{D29.943}
Since we're after an understanding of the functor $f^\times$, we
need to agree what the word ``understanding'' will mean. Recall that  
the adjunction $\R f_*\dashv f^\times$ gives a natural isomorphism
\[\xymatrix@C+30pt{
\Hom(A,f^\times B)\ar[r]^-{\ph(A,B)} & \Hom(\R f_*A,B)
}\]
Putting $A=f^\times B$ this specializes to a homomorphism
\[
\xymatrix@C+30pt{
\Hom(f^\times B,f^\times B)\ar[r]^-{\ph(f^\times B,B)} & \Hom(\R f_*f^\times B,B)
}\]
which sends $\id:f^\times B\la f^\times B$ to the map 
$\e:\R f_*f^\times B\la B$, the counit of adjunction. Consider 
\[
\xymatrix{
\Hom(A,f^\times B)\ar[r]^-{\R f_*} & \Hom(\R f_*A,\R f_*f^\times
B)\ar[rr]^-{\Hom(-,\e)} & &\Hom(\R f_*A,B)
}\]
It is classical that naturality forces the composite to agree with
$\ph(A,B)$. Summarizing:
\edis

\cnc{C0.1}
If we could compute, for every $B\in\Dqc(Y)$, the object $f^\times B$ and
the morphism $\e:\R f_*f^\times B\la B$, then we'd feel we understand
the adjunction pretty well. After all the map
$\ph(A,B):\Hom(A,f^\times B)\la\Hom(\R f_*A,B)$ would be explicit: given
an element
$\alpha\in\Hom(A,f^\times B)$, that is a morphism $\alpha:A\la f^\times B$,
then the map $\ph(A,B)$ would send $\alpha$ to the
composite $\R f_*A\stackrel{\R f_*\alpha}\la \R f_*f^\times B\stackrel\e\la B$,
which is an element $\e\circ\R f_*\alpha\in\Hom(\R f_*A,B)$.
OK: it wouldn't be so clear how to go back, but classically people have
been happy with understanding just this direction.
\ecnc

We will soon specialize to the case where $f$ is smooth and proper,
but the next result holds more generally and we state it in a strong form.

\thm{T0.3}
Assume $f:X\la Y$ is a finite-type morphism of noetherian schemes.
If $B\in\Dqc(Y)$ and $C\in\Dqc(X)$ then there is a canonical
natural isomorphism
$p_{B,C}^{}:B\oo^\LL \R f_*C\la
  \R f_*(\LL f^*B\oo^\LL C)$
and  a canonical natural transformation
$\chi:\LL f^*B\oo^\LL f^\times\co_Y^{}\la f^\times B$
such that the following pentagon commutes
\[
\xymatrix@C-20pt{
\R f_*\big[\LL f^*B\oo^\LL f^\times\co_Y^{}\big]\ar[dr]_-{\R f_*\chi}\ar[rr]^-{p_{B,f^\times\co_Y^{}}^{-1}}& & B  \oo_{}^\LL\R f_*f^\times \co_Y^{} \ar[rr]^-{\id\oo\e} &&  B\oo^\LL_{}\co_Y^{}\ar@{=}[dl]\\
& \R f_*f^\times B\ar[rr]_-\e && B &
}\]
Furthermore: the map $\chi$ is an isomorphism if and only if $f$ is proper
and of finite
Tor-dimension.
\ethm

The non-expert should view finite Tor-dimension as a technical condition that
will be satisfied by all the $f$'s we will consider. We will discuss
properness in Remark~\ref{R9999.1010}.

\rmk{R-1.5}
Suppose  $f$ is proper and of finite Tor-dimension.
Then we have an isomorphism $\LL f^*B\oo_{}^\LL f^\times\co_Y^{}\la f^\times B$,
and the commutative pentagon
of Theorem~\ref{T0.3}
makes precise
the compatibility of this
isomorphism with the counit $\e$ of the adjunction $\R f_*\dashv f^\times$.
Thus the import of Theorem~\ref{T0.3} is that,
as long as $f$ is proper and of finite Tor-dimension,
it suffices to compute
$f^\times\co_Y^{}$ and the counit of adjunction
$\e:\R f_*f^\times\co_Y^{}\la\co_Y^{}$. We are reduced to studying a single
object, namely $\co_Y^{}\in\Dqc(Y)$.
\ermk

\rmd{R-1.7}
The next reduction comes from the observation that the category
$\Dqc(X)$ has many endofunctors. 
There are  many diagrams
\[\xymatrix@C+20pt@R-20pt{
   & \ar@{=>}[dd]_{c_W^{}} & \\
  \Dqc(X)\ar@/^1.5pc/[rr]^-{\Gamma_W^{}}\ar@/_1.5pc/[rr]_-\id & &\Dqc(X)\\
  & &
}\]
That is: there are many choices of functors $\Gamma_W^{}:\Dqc(X)\la\Dqc(X)$,
which come together with natural transformations $c_W^{}:\Gamma_W^{}\la\id$.
The ones we have in mind are the Bousfield colocalizations. They
come about as follows.

For every point $p\in X$ let $i_p:p\la X$ be the inclusion,
which we view as a morphism of schemes
$i_p:\text{\rm Spec}\big(k(p)\big)\la X$. 
Suppose we are given a set of points $W\subset X$.
The full subcategory $\D_{\mathbf{qc},W}^{}(X)\subset\Dqc(X)$ will be the
subcategory of all objects \emph{supported on $W$,} we recall
that this means
\[
\D_{\mathbf{qc},W}^{}(X)\eq\{E\in\Dqc(X)\mid \LL i_p^*E=0\text{ for all }p\notin W\}\ .
\]
Let $I_W:\D_{\mathbf{qc},W}^{}(X)\la\Dqc(X)$ be the inclusion.
A straightforward generalization of a theorem of
Bousfield~\cite{Bousfield79A} tells us that $I_W$ has a right adjoint
$J_W:\Dqc(X)\la\D_{\mathbf{qc},W}^{}(X)$, and the colocalizations we have in mind are
the counits of adjunction $c_W^{}:I_WJ_W\la\id$.

We will give a concrete example later in this section.
\ermd

\rmk{R-1.999}
Let the notation be as in Reminder~\ref{R-1.7}. It is customary to choose
the set of points $W\subset X$ to be \emph{closed under specialization},
meaning if $p\in W$, and if $q\in X$ belongs to the closure $\ov{\{p\}}$ of
$p$, then $q\in W$. The advantage is that for such $W$ the Bousfield
colocalization can be computed locally. More precisely: let $u:U\la X$ be an
open immersion, let $\D_{\mathbf{qc},U\cap W}^{}(U)\subset\Dqc(U)$ be
the full subcategory of
objects supported on $U\cap W$, let
$I_{U\cap W}^{}:\D_{\mathbf{qc},U\cap W}^{}(U)\la\Dqc(U)$
be the inclusion and $J_{U\cap W}^{}$ its right adjoint,
and let $c_{U\cap W}^{}:\Gamma_{U\cap W}^{}\la \id_{\Dqc(U)}^{}$ be the counit of
the adjunction $I_{U\cap W}^{}\dashv J_{U\cap W}^{}$. Then the
relation with the $c_W^{}:\Gamma_W^{}\la\id_{\Dqc(X)}^{}$ of
Reminder~\ref{R-1.7} is simple:
there is a
canonical isomorphism $\LL u^*\Gamma_W^{}\cong\Gamma_{U\cap W}^{}\LL u^*$ making
the triangle below commute
\[\xymatrix@R-20pt@C+40pt{
\LL u^*\Gamma_W^{} \ar[rd]^-{\LL u^*c_W^{}}\ar[dd]_-{\wr}& \\
  & \LL u^*\\
\Gamma_{U\cap W}^{}\LL u^*\ar[ru]_-{c_{U\cap W}^{}\LL u^*}
}\]
If $W$ isn't specialization-closed this may fail, 
and the bottom line is that we feel infinitely
more comfortable working with tools
that lend themselves to local computations.
\ermk

Given a $c_W^{}:\Gamma_W^{}\la\id$ as in Reminder~\ref{R-1.7},
we can form the next gadget:

\dfn{D-1.9}
We define $\rho_W^{}$ to be the composite
\[
\xymatrix@C+20pt{
\R f_*\Gamma_W^{} f^\times\co_Y^{}\ar[rr]^-{\R f_*c_W^{}f^\times}
& &\R f_*f^\times\co_Y^{}\ar[r]^-\e &\co_Y^{}
}\]
\edfn

\rmk{R-1.11}
Ideally, we would aim to choose $c_W^{}:\Gamma_W^{}\la\id$ in such a way that
\be
\item
  The composite $\rho_W^{}$ of Definition~\ref{D-1.9} is easy to compute.
\item
  From the computation of $\rho_W^{}$ we learn a lot about
  $\e:\R f_*f^\times\co_Y^{}\la\co_Y^{}$.
\ee
Of course we could make a dumb choice of $c_W^{}:\Gamma_W^{}\la\id$. For
example: if we let $c_W^{}:\Gamma_W^{}\la\id$
be the identity map $\id\la\id$, then $\rho_W^{}=\e$, 
we don't lose any information in passing from $\e$ to  $\rho_W^{}$,
but we also haven't simplified the computation.
Or if we choose $\Gamma_W^{}=0$
then the computation of $\rho_W^{}$ becomes trivial,
but worthless. The important thing is to choose $c_W^{}:\Gamma_W^{}\la\id$
wisely.
\ermk

\subsection{If $f$ is smooth and proper}
\label{S3}
In the most classical case of the theory we have the following results:

\thm{T0.13}
Assume $f:X\la Y$ is smooth and proper, of relative dimension $n$.
Then there is a canonical isomorphism
$\theta:\Omega_{f}^n[n]\stackrel\sim\la f^\times\co_Y^{}$.
\ethm

\rmk{R9999.1010}
We should explain the theorem, starting with the
hypotheses: if the non-expert tried to guess what it means for $f$ to be
smooth and proper, chances are she was right about proper but wrong about
smooth. Let us elaborate.

It is customary to consider the following two conditions, which a
continuous map $f:X\la Y$  of topological spaces can satisfy:
\be
\item
$f^{-1}(K)$ is compact whenever $K\subset Y$ is compact.
\item
The map $f$ is universally closed. This means that, if $f':X'\la Y'$
is some pullback of $f$, then $f'(K)\subset Y'$ is closed whenever
$K\subset X'$ is closed.
\ee
In the category of locally compact Hausdorff spaces the two
are equivalent, and a map satisfying these equivalent conditions
is what's normally called proper.
As it happens the
topological spaces that come up in algebraic geometry are rarely
Hausdorff, and in the category of schemes (i) and (ii) aren't
equivalent. It turns out that the right way to define proper
maps of schemes is to use (ii), this yields the theory one
would intuitively expect.

But when it comes to smoothness algebraic geometers
chose to be contrary. In differential geometry---and hence also
in related topics like PDE---a smooth map of manifolds is
defined to be a $C^\infty$
map. With this definition alegbraic geometers never consider
any map that's remotely non-smooth.

Even though the term ``smooth map'' was already in use
in a well-established, 
clearly delineated context, algebraic geometers decided to steal 
the word
and give it a different meaning. In this survey we follow
the conventions of algebraic geometry:
what we
label a ``smooth map'' is what everyone else would dub a ``submersion''. In
algebraic geometry, a morphism $f:X\la Y$ of manifolds is \emph{smooth}
if, at every point $p\in X$, the derivative is a surjection
$t_p:T_p\la T_{f(p)}$. Here $T_p$ is the tangent space
at $p$ and $T_{f(p)}$ is the tangent space at $f(p)$. The smooth
map $f$ has \emph{relative dimension $n$} if the
kernel of the linear map $t_p:T_p\la T_{f(p)}$ is $n$--dimensional for
every $p\in X$.

Assume $f$ is a smooth map of
relative dimension $n$ in the sense above, and let $\Omega_X^{},\Omega_Y^{}$
be
the cotangent bundles of $X,Y$ respectively.
The pullback $f^*\Omega_Y^{}$ is naturally a subbundle
of $\Omega_X^{}$, and the \emph{relative cotangent bundle} is by
definition the quotient
$\Omega_{f}^{}=\Omega_X^{}/f^*\Omega_Y^{}$, often
written $\Omega_{f}^{1}$. It is a 
vector bundle of  rank $n$ over $X$, whose top exterior power is what is
usually written $\Omega_{f}^{n}=\wedge^n\Omega_{f}^{1}$. The line bundle
$\Omega_{f}^{n}$ is called the \emph{relative canonical bundle} of $f$.
Thus Theorem~\ref{T0.13}
is the assertion that $f^\times\co_Y^{}$ is canonically isomorphic to
the object $\Omega_{f}^{n}[n]\in\Dqc(X)$, which is 
the cochain
complex with only one nonvanishing term, namely
the relative canonical bundle $\Omega_{f}^n$ in degree $-n$.
\ermk

\rmk{R0.15}
With the notation as explained in Remark~\ref{R9999.1010},
Theorem~\ref{T0.13}
computes for us the object $f^\times\co_Y^{}$. By Remark~\ref{R-1.7}
our mission would be accomplished if we could also compute
the counit of adjunction $\e:\R f_*f^\times\co_Y^{}\la\co_Y^{}$.
In Remark~\ref{R-1.11} we noted that it might prove expedient to
take advantage of some Bousfield colocalization $c_W^{}:\Gamma_W^{}\la\id$.
The traditional choice, which happens to be
well-suited for the current computation, is to take
$c_W^{}:\Gamma_W^{}\la \id$ to be the Bousfield
colocalization of Reminder~\ref{R-1.7}, where the set of points $W\subset X$
is the union of the irreducible closed subsets
$Z\subset X$
such that the composite map $Z\la X\la Y$ is generically finite.
\ermk

It should be noted that our $\Gamma_W^{}$ has been extensively studied and
is very computable, the subject dealing with functors of this genre is called
\emph{local cohomology.} Most of what's written about $\Gamma_W^{}$ is in the
commutative algebra literature.

Before the theorem it might help to
illustrate the abstraction in a simple case.

\exm{E0.17}
Suppose $Y=\text{Spec}(k)$ where $k$ is field.
Then $X$ is smooth over the field $k$ and $n$--dimensional. Take
a minimal injective resolution for $\Omega_{f}^n[n]\cong f^\times\co_Y^{}$:
that is
form a cochain complex
\[\xymatrix{
  0\ar[r] & I^{-n}\ar[r] &I^{-n+1}\ar[r] &
  \cdots\ar[r] & I^{-1}\ar[r] & I^0\ar[r] & 0
}\]
where $I^{-n}$ is the injective
envelope of $\Omega_{f}^n$, next $I^{-n+1}$ is the injective envelope
of $I^{-n}/\Omega_{f}^n$, and so on.\footnote{
The reader might recall that injective resolutions, even minimal ones, are
only unique up to homotopy---hence it isn't obvious that the zeroth 
sheaf $I^0$ 
is functorial in anything---and the formula that's about to come, that is
$I^0=\Gamma_W^{}\Omega_{f}^n[n]$, seems unreasonable at first sight.
 
Given any two injective resolutions
of $I^*$ and $J^*$ of $\Omega_{f}^n[n]$ there are cochain 
maps
$I^*\stackrel\alpha\la J^*\stackrel\beta\la I^*$, unique
up to homotopy, so that $\alpha\beta$ and
$\beta\alpha$ are homotopic to the identity. What is special here,
because we are dealing with minimal injective resolutions of a line
bundle 
$\Omega^n_{f}$ on a regular scheme $X$, is that any homotopy must vanish.
Hence the minimal injective resolution $I^*$ is unique up to canonical
isomorphism, as is $I^0=\Gamma_W^{}\Omega_{f}^n[n]$.

For the experts in commutative algebra: the way one proves the vanishing of
any homotopy is by noting that the 
injective sheaf $I^{-j}$ is a direct sum of indecomposable injectives
supported at points of dimension $j$, and there are no non-zero maps
from an indecomposable injective supported at a point of dimension
$j$ to an indecomposable injective supported at a point of dimension $j+1$.}
 Note that $X$ is regular and
$n$--dimensional, hence the injective dimension of $\Omega_{f}^n$
is $\leq n$. The minimal
injective resolution must stop no later than $I^0$.
The $W$ of Remark~\ref{R0.15} is the set of all closed points
in $X$.
The corresponding Bousfield colocalization
$c_W^{}:\Gamma_W^{}\Omega_{f}^n[n]\la\Omega_{f}^n[n]$
comes down to
the cochain map
\[\xymatrix{
   0\ar[r] & 0\ar[r]\ar[d] &0\ar[r]\ar[d]  &
  \cdots\ar[r] & 0\ar[d] \ar[r] & I^0\ar@{=}[d] \ar[r] & 0\\
  0\ar[r] & I^{-n}\ar[r] &I^{-n+1}\ar[r] &
  \cdots\ar[r] & I^{-1}\ar[r] & I^0\ar[r] & 0
}\]
In the next theorem we will compute the composite
$\rho_W^{}$ of Definition~\ref{D-1.9}, that is we propose to compute the map
$\R f_*\Gamma_W^{}f^\times\co_Y^{}\stackrel{\R f_*c_W^{}}\la
\R f_*f^\times\co_Y^{}\stackrel\e\la\co_Y^{}$.
This map identifies, via
the isomorphism $\Omega^n_{f}[n]\cong f^\times\co_Y^{}$ of
Theorem~\ref{T0.13}, with a composite
$\R f_*\Gamma_W^{}\Omega^n_{f}[n]\stackrel{\R f_*c_W^{}}\la
\R f_*\Omega^n_{f}[n]\stackrel\e\la\co_Y^{}$.
Note that as $Y=\spec k$ we have $\co_Y^{}=k$.
Recalling
Reminder~\ref{R-1.-1}(ii), we compute $\R f_*$ by applying $f_*$ to an injective
resolution. That is, in our case $\rho_W^{}$ comes down to a composite
cochain map
\[\xymatrix{
   0\ar[r] & 0\ar[r]\ar[d] &0\ar[r]\ar[d]  &
  \cdots\ar[r] & 0\ar[d] \ar[r] & f_*I^0\ar@{=}[d] \ar[r] & 0\\
  0\ar[r] & f_*I^{-n}\ar[r]\ar[d] &f_*I^{-n+1}\ar[r]\ar[d] &
  \cdots\ar[r] & f_*I^{-1}\ar[r]\ar[d] & f_*I^0\ar[d]^{\e^0}\ar[r] & 0\\
   0\ar[r] & 0\ar[r] &0\ar[r]  &
  \cdots\ar[r] & 0 \ar[r] & k \ar[r] & 0
}\]
This makes it obvious why the
$c_W^{}:\Gamma_W^{}f^\times\co_Y^{}\la f^\times\co_Y{}$ of Remark~\ref{R0.15}
is a reasonable choice: it certainly doesn't lose information,
the map $\e^0:f_*I^0\la k$ most definitely
determines the cochain map $\e:\R f_*\Omega^n_{f}[n]=f_*I^*\la k$.

Injective resolutions might be useful for proving a map is 
informative,
but for computations one usually prefers other resolutions.
Fortunately the functor $\Gamma_W^{}$ has other descriptions.
For example there is a description in terms of local cohomology:
it turns out that the
$f_*I^0=\R f_*\Gamma_W^{}\Omega_{f}^{}[n]$
above is isomorphic to
\[
f_*I^0\quad\cong\quad\bigoplus_{p\in X,\,\, \ov{\{p\}}=\{p\}}H^n_p(\Omega^n_{f})\ .
\]
That is, $f_*I^0=\R f_*\Gamma_W^{}\Omega_{f}^n[n]$ is
the direct sum, over all closed points $p\in X$, of the $n\mth$ local
cohomology of the sheaf $\Omega_{f}^n$ at $p$. The standard
computation of local cohomology, via \u{C}ech complexes, tells
us that $f_*I^0$ may
be written as the quotient of $J$, where $J$ is the direct
sum, over all closed points $p\in X$, of the vector space of
meromorphic differential forms at $p$. A \emph{meromorphic form
  at $p$} is the following list of data:
\be
\item The closed point $p\in X$. 
\setcounter{enumiv}{\value{enumi}}
\ee
For parts (ii) and (iii) below we write
$\co^{}_{X,p}$ for the stalk at $p$ of the structure sheaf
$\co_X^{}$.
\be
\setcounter{enumi}{\value{enumiv}}
\item
  A system of coordinates $(x_1^{},\ldots,x_n^{})$ at $p$, that
  is $(x_1^{},\ldots,x_n^{})$ generate the maximal ideal of the
  local ring $\co^{}_{X,p}$.
\item
  An expression
  \[\frac{fdx_1^{}\wedge dx_2^{}\wedge\cdots\wedge dx_n^{}}{x_1^{N}x_2^{N}\cdots x_n^{N}}\]
  where $f\in\co^{}_{X,p}$.
\ee
\eexm

In Example~\ref{E0.17} we simplified our life by assuming $Y=\spec k$ with
$k$ a field. The general case is slightly
more cumbersome to describe but similar. And now for the main result.

\thm{T0.19}
Suppose $f:X\la Y$ is a smooth and proper morphism of noetherian
schemes, identify $\Omega^n_{f}[n]\cong f^\times\co^{}_Y$
via the canonical isomorphism $\theta$ of Theorem~\ref{T0.13},
and let
$c_W^{}:\Gamma_W^{}\la\id$ be as chosen in
Remark~\ref{R0.15}. Then the map
$\rho_W^{}:\R f_*\Gamma_W^{}\Omega^n_{f}[n]\la\co_Y^{}$
of Definition~\ref{D-1.9} may be represented, in the
notation of Example~\ref{E0.17} [more accurately in the generalization
  of the notation to the case where $Y$ is arbitrary],
by the map $J\la \co_Y^{}$ taking a finite sum of
meromorphic forms
to the sum of their residues.
\ethm

\exm{E0.21}
Let us return to the situation of Example~\ref{E0.17}, where $Y=\spec k$.
But now assume further that
$k=\ov k$ is
an algebraically closed field. Then the residue of a meromorphic
differential form is the obvious. In the expression
\[\frac{fdx_1^{}\wedge dx_2^{}\wedge\cdots\wedge dx_n^{}}{x_1^{N}x_2^{N}\cdots x_n^{N}}\]
we may expand $f\in\co^{}_{X,p}$ into a Taylor series, which we view
as an element of
the completion $\widehat\co^{}_{X,p}$
of the ring $\co^{}_{X,p}$. This gives
an expansion of the entire meromorphic form into a Laurent series 
\[\sum_{\text{all }k_i\leq N}\frac{a_{k_1,k_2,\ldots,k_n}^{}dx_1^{}\wedge dx_2^{}\wedge\cdots\wedge dx_n^{}}{x_1^{k_1}x_2^{k_2}\cdots x_n^{k_n}}\]
where the coefficients $a_{k_1,k_2,\ldots,k_n}^{}$ belong to the
field $k$. The map $\rho_W^{}$ takes the meromorphic form to $a_{1,1,\ldots,1}^{}$,
that is to the coefficient of
\[\frac{dx_1^{}\wedge dx_2^{}\wedge\cdots\wedge dx_n^{}}{x_1^{}x_2^{}\cdots x_n^{}}\]
in the Laurent series.
\eexm

\rmk{R0.23}
The discussion of
Example~\ref{E0.21} should explain why we made the simplifying assumption
that $Y=\spec k$ with $k$ algebraically closed. We run into subtleties
already when
$Y=\spec k$ but we drop the hypothesis that $k=\ov k$.
In this case the
closed point $p$ will not in general be $k$-rational,
and we would not expect an element
$f\in\co^{}_{X,p}\subset\widehat\co^{}_{X,p}$ to have a Taylor expansion
with coefficients in $k$,
in the generators $\{x_1^{},x_2^{},\ldots ,x_n^{}\}$ of
the maximal ideal. The definition of
the residue of a meromorphic form becomes subtler.
\ermk

\subsection{Application: Serre duality}
\label{SS29.987}

Serre's classical duality theorem is a special case; let us recall precisely
how. But first a reminder: back in Remark~\ref{R24.989796}
we disclosed that the non-expert will be asked to accept, on faith, all
computations of Hom-sets in derived categories---for example the ones
she's about to witness.

Continue to assume that $k$ is a field and
$f:X\la Y=\spec k$ is smooth and proper, of relative dimension $n$.
Let $\cv$ be a vector
bundle on $X$.
Then the
adjunction $\R f_*\dashv f^\times$ tells us that, in the commutative square
below, the map $\ph$ is an isomorphism
\[\xymatrix@C+30pt{
   \Hom_{\Dqc(X)}^{}\big(\cv[i],f^\times\co_Y^{}\big)\ar@{=}[d]
   \ar[r]^-\ph& \Hom_{\Dqc(Y)}^{}\big(\R f_*\cv[i],\co_Y^{}\big) \ar@{=}[d]\\
\Hom_{\Dqc(X)}^{}\big(\cv[i],\Omega^n_X[n]\big)\ar[r] & \Hom_{\Dqc(Y)}^{}\big(\R f_*\cv[i],k\big)
}\]
The bottom row simplifies to
\[
H^{n-i}\big(\HHom(\cv,\Omega_{X}^n)\big)\quad\cong \quad
\Hom\big(H^i(\cv),k\big)
\]
meaning there is a nondegenerate pairing
$H^i(\cv)\oo H^{n-i}\big(\HHom(\cv,\Omega_{X}^n)\big)\la k$.
The pairing is explicit.
It takes a morphism $\alpha:\co_X^{}\la\cv[i]$ and a
morphism $\beta:\cv\la\Omega_X^n[n-i]$ to the image of
$1\in\R f_*\co_X^{}$ under the
composite
\[\xymatrix@C+20pt{
\R f_*\co_X^{}\ar[r]^-{\R f_*\alpha} &
\R f_*\cv[i] \ar[r]^-{\R f_*\beta[i]} &
\R f_*\Omega_X^n[n]\ar[r]^-{\e} & k 
}\]
After all: Discussion~\ref{D29.943} and Conclusion~\ref{C0.1}
combine to tell us 
that the map
$\beta[i]:\cv[i]\la\Omega_X^n[n]\cong f^\times\co_Y^{}$ goes under 
the bijection $\ph$ to
the composite
\[\xymatrix@C+20pt{
\R f_*\cv[i] \ar[r]^-{\R f_*\beta[i]} &
\R f_*\Omega_X^n[n]\ar[r]^-{\e} & k 
}\]
and all we do is evaluate
this composite at the element $\alpha:\co_X^{}\la\cv[i]$ of the vector space
$H^i(\cv)=H^i\big(\HHom(\co_X^{},\cv)\big)=\Hom\big(\co_X^{},\cv[i]\big)$.

\rmk{R990.999}
We can view Grothendieck duality as being Serre duality on steroids.
Being macho, Grothendieck duality doesn't restrict the scheme
$Y$ to be the one-point space,
doesn't assume the map $f$ to be smooth and proper, and doesn't confine itself
to only dealing with
vector bundles.
\ermk

\section{The proofs}
\label{S30}

Modulo technicalities, the theorems of Section~\ref{S29} are all contained
in Hartshorne~\cite{Hartshorne66}. More precisely: the theorems in
Hartshorne~\cite{Hartshorne66} aren't quite as clean or general,
Section~\ref{S29} is comprised of several technical improvements on
those original assertions.
However: with one inessential exception---which we
will mention at the very end of Sketch~\ref{SK30.3}---all the improvements
had been obtained by the mid-1990s. In other words:
none of the results in
Section~\ref{S29} is younger than two decades.

What is new is that we can now prove \emph{every one of these statements}
simply and
directly, sidestepping the customary circuitous routes and
long detours, and bypassing all of the traditional stopovers on distant
planets.
We will next discuss where the reader can find these simple, formal proofs.
This naturally divides into two parts.

\subsection{Simple proofs that have been around for decades}
\label{SS30.1}

Let us begin with Reminder~\ref{R-1.-1}: we gave an explicit construction
of the functors $\LL f^*$ and $\R f_*$, and asserted the existence of
a functor $f^\times$ right adjoint to $\R f_*$. The first short and formal
proof of the existence of $f^\times$ may be found in Deligne's
appendix~\cite{Deligne66}
to Hartshorne's
book~\cite{Hartshorne66}. In that proof the schemes are assumed
noetherian and the derived categories are of bounded below complexes.
The reader may find more general theorems in
Balmer, Dell'Ambrogio and Sanders~\cite{Balmer-DellAmbrogio-Sanders15}
and in \cite{Neeman96,Neeman13},
with the strongest theorem to date covering the case where $f$ is any
concentrated morphism of quasicompact, quasiseparated algebraic stacks.
The modern proofs work by showing that $\R f_*$
respects coproducts and applying Brown representability.

Now we turn to Theorem~\ref{T0.3}, and for the reader's convenience we recall
the statement

\bigskip

\nin
{\bf Theorem~\ref{T0.3}, [reminder].}\ \ \emph{
Assume $f:X\la Y$ is a finite-type morphism of noetherian schemes.
Then there is a canonical
natural isomorphism
$p_{A,B'}^{}:A\oo^\LL \R f_*B'\la
  \R f_*(\LL f^*A\oo^\LL B')$
and  a canonical natural transformation
$\chi:\LL f^*A\oo^\LL f^\times\co_Y^{}\la f^\times A$
such that the following pentagon commutes
\[
\xymatrix@C-20pt{
\R f_*\big[\LL f^*A\oo^\LL f^\times\co_Y^{}\big]\ar[dr]_-{\R f_*\chi}\ar[rr]^-{p_{A,f^\times\co_Y^{}}^{-1}}& & A  \oo_{}^\LL\R f_*f^\times \co_Y^{} \ar[rr]^-{\id\oo\e} &&  A\oo^\LL_{}\co_Y^{}\ar@{=}[dl]\\
& \R f_*f^\times A\ar[rr]_-\e && A &
}\]
Furthermore: the map $\chi$ is an isomorphism if and only if $f$ is proper
and of finite
Tor-dimension.}

\bigskip

\nin
The modern proof may be found
in \cite{Neeman96}, and the reader might also wish to look at
\cite{Balmer-DellAmbrogio-Sanders15} for a generalization of
Theorem~\ref{T0.3} proved by the same techniques.
To emphasize the formal nature of the argument we give an outline.

\skt{SK30.3}
The functor $\LL f^*$ is strong monoidal---it respects the tensor
product, there is a natural isomorphism
$\LL f^*(A\oo^\LL B)\la (\LL f^*A)\oo^\LL (\LL f^*B)$.
If we put $B=\R f_*B'$ this gives the first map in
the composite
\[\xymatrix{
\LL f^*(A\oo^\LL \R f_*B')\ar[r] & (\LL f^*A)\oo^\LL (\LL f^*\R f_*B')
\ar[rr]^-{\id\oo\e'} && (\LL f^*A)\oo^\LL B'\ ,
}\]
where $\e':\LL f^*\R f_*\la\id$ is the counit of the adjunction
$\LL f^*\dashv\R f_*$. Adjunction, applied to
the highlighted composite above, gives a corresponding
map $p_{A,B'}^{}:A\oo^\LL \R f_*B'\la \R f_*\big[(\LL f^*A)\oo^\LL B'\big]$.
The map $p$ is an isomorphism, the so-called
classical ``projection formula''. But since we are into seeing
what part of the theory is formal, let us indicate the modern proof
that $p$ is an isomorphism.

If $A$ is a perfect complex then it is
``strongly dualizing''\footnote{Recall that an object $A$
in a monoidal category  is strongly dualizing
if there exists a dual object $A^\vee$ and maps $A^\vee\oo A\la\one$ and
$\one\la A\oo A^\vee$ so that the composites $A\la A\oo A^\vee\oo A\la A$
and $A^\vee\la A^\vee\oo A\la A^\vee$ are both the identity.
For the monoidal category $\Dqc(X)$ the strongly dualizing objects
are the perfect complexes.}---in particular there
exists a dual complex $A^\vee$ and a canonical isomorphism
$A\oo^\LL_{}(-)\cong\RHHom(A^\vee,-)$.
The following string of isomorphisms
\begin{eqnarray*}
\Hom(C\,\,,\,\,A\oo^\LL \R f_*B') &\cong& \Hom\big(C\,\,,\,\,\RHHom(A^\vee, \R f_*B')\big)  \\
&\cong&\Hom(C\oo^\LL A^\vee\,\,,\,\, \R f_*B')\\
&\cong&\Hom\big((\LL f^*(C\oo^\LL A^\vee)\,\,,\,\, B'\big)\\
&\cong&\Hom\big((\LL f^*C)\oo^\LL(\LL f^* A^\vee)\,\,,\,\, B'\big)\\
&\cong&\Hom\big((\LL f^*C)\oo^\LL(\LL f^* A)^\vee\,\,,\,\, B'\big)\\
&\cong&\Hom\Big[(\LL f^*C)\,\,,\,\, \RHHom\big((\LL f^* A)^\vee,B'\big)\Big]\\
&\cong&\Hom\Big[C\,\,,\,\, \R f_*\big((\LL f^* A)\oo^\LL B'\big)\Big]
\end{eqnarray*}
holds for every $C$ and is a formal consequence of the definition of
strongly dualizable objects---in particular the fifth
isomorphism is formal, any strict monoidal functor must
take strongly dualizable objects to strongly
dualizable objects, and must respect duals. Yoneda gives that the isomorphism
of Hom-sets must come from an isomorphism
$A\oo^\LL \R f_*B'\cong \R f_*\big(\LL f^* A)\oo^\LL B'\big)$, and it is
an exercise in the definitions to check that this isomorphism is induced
by the map $p_{A,B'}^{}$. In other words: the map
$p_{A,B'}^{}$ induces an isomorphism as
long as $A$ is strongly dualizing---in the category
$\Dqc(Y)$ this means as long as $A$ is a perfect complex.

Thus the subcategory of all $A$'s for which the map $p_{A,B'}^{}$
induces an isomorphism, for every $B'$, contains
the perfect complexes and is closed under suspensions, triangles and
coproducts. The closure under triangles and coproducts is because the
functors $\LL f^*$
and $\R f_*$ both respect triangles and
coproducts. But then \cite[Lemma~3.2]{Neeman96}
tells us that
every object $A\in\Dqc(Y)$ belongs---the map $p_{A,B'}^{}$ is an isomorphism
for all pairs $A,B'$.

Once we know that the map $p$ is an isomorphism we can repeat the
idea. Put $B'=f^\times B''$, and then $p^{-1}_{A,f^\times B''}$ gives the first map
in the composite
\[\xymatrix{
\R f_*\big[(\LL f^*A)\oo^\LL (f^\times B'')\big]\ar[rr]^-{p^{-1}_{A,f^\times B''}} &&
A\oo^\LL \R f_*f^\times B''\ar[rr]^-{\id\oo\e} & & A\oo^\LL B''
}\]
In this composite
$\e:\R f_*f^\times B''\la B''$ is the counit of the adjunction
$\R f_*\dashv f^\times$. Adjunction tells us that the composite corresponds
to a map
$\chi(A,B''):(\LL f^*A)\oo^\LL (f^\times B'')\la f^\times(A\oo^\LL B'')$.
The map $\chi=\chi_A^{}$ of Theorem~\ref{T0.3} is just the special
case where $B''=\co_Y^{}$. The commutativity of the pentagon of
Theorem~\ref{T0.3} simply spells out what it means for the map $\chi_A^{}$
to correspond, under the adjunction $\R f_*\dashv f^\times$, to the
composite above---in other words the pentagon commutes by definition.

It remains to discuss when the map $\chi$ is an isomorphism.

One can easily write down a string of isomorphisms, much like those
we used above to study $p_{A,B'}^{}$,
which combine to show that the map $\chi_A^{}$ is an isomorphism as
long as $A$ is a perfect complex; see \cite[top of page~228]{Neeman96}.
Consider the category of all $A$'s so that the map $\chi_A^{}$
is an isomorphism---it contains the perfect complexes, is closed under
suspensions and triangles, and also closed under coproducts
if $f^\times$
respects coproducts. Thus, as
long as $f^\times$ respects coproducts, the map $\chi_A^{}$ is an isomorphism
for every $A$. In fact
the condition that $f^\times$ respects coproducts is necessary and
sufficient for the map $\chi$ to be an isomorphism.

Now \cite[Lemma~5.1]{Neeman96} comes to our aid:
given a pair of adjoint triangulated functors
$F\dashv G$, between compactly generated triangulated categories, the right
adjoint $G$ respects coproducts if and only if the left adjoint $F$
respects compact
objects.\footnote{Let $\ct$ be a triangulated category with coproducts.
  An object $C\in\ct$ is \emph{compact} if $\Hom(C,-)$ respects coproducts---for
  this survey
  what's relevant is that the compact objects in $\Dqc(X)$ are the perfect
  complexes.}
Specializing to the pair of adjoint functors $\R f_*\dashv f^\times$, this
tidbit of formal nonsense says that
the functor $f^\times$ respects coproducts
if and only if $\R f_*$ takes perfect complexes to perfect complexes.
Thus we have tranformed a question about the
right adjoint $f^\times$, which is mysterious, into a problem about
its left adjoint $\R f_*$, which is explicit and computable. It is an
old theorem of Illusie~\cite[Expos\'e~III, Corollaire~4.3.1(a)]{Illusie71C}
that if $f$ is proper and of finite Tor-dimension then $\R f_*$ respects
perfect complexes. The converse [which we don't use in this
article], that is the theorem that if $\R f_*$ respects 
perfect complexes then $f$ must be proper and of finite Tor-dimension,
is much more recent. It may be found in \cite{Lipman-Neeman07}. 
\eskt

This concludes our discussion of~Theorem~\ref{T0.3}.
Before proceeding to the remaining two theorems we include a brief
interlude, recalling
a couple of base-change maps. 

\rmd{R30.30303}
Let us begin at the formal level. Suppose we are given four categories
$\cw$, $\cx$, $\cy$ and $\cz$, as well as 
pairs of adjoint functors
\[
\xymatrix{\gamma:\cy\ar@<0.5ex>[rr]& &\cw:\Gamma\ar@<0.5ex>[ll] &\text{ and } &
\beta:\cz\ar@<0.5ex>[rr]& &\cx:B\ar@<0.5ex>[ll]}
\]
Recall: this is shorthand for the
assertion that $\Gamma:\cw\la\cy$ is right adjoint to $\gamma:\cy\la\cw$,
and $B:\cx\la\cz$ is right adjoint to $\beta:\cz\la\cx$.
Assume further that we are given a pair of functors
\[
\xymatrix{\alpha:\cx\ar[rr]& &\cw &\text{ and } &
\delta:\cz\ar[rr]& &\cy}
\]
Then there is a canonical bijection between natural transformations
\[\xymatrix@R-22pt@C+20pt{
  \cw&\ar[l]_-{\alpha} \cx & & \cw\ar[dddd]_-{\Gamma}&\ar[l]_-{\alpha} \cx \ar[dddd]^-{B}\\
\quad  &\quad & &\quad &\quad \\
\quad&\quad&\text{and} &\quad &\quad \\
\quad \ar@{=>}[ruu]^\rho &\quad& & \quad &\quad \ar@{=>}[uul]_\sigma  \\
 \cy\ar[uuuu]^-{\gamma} &\ar[l]^-{\delta}  \cz\ar[uuuu]_-{\beta}& & \cy &
 \cz\ar[l]^-{\delta} 
}\]
The formula is explicit: if $\eta(\gamma\dashv\Gamma)$ and
$\e(\gamma\dashv\Gamma)$ are, respectively, the unit and counit of the
adjunction $\gamma\dashv\Gamma$, while
$\eta(\beta\dashv B)$ and $\e(\beta\dashv B)$
are, respectively, the unit and counit of the
adjunction $\beta\dashv B$, then the formulas are that the bijection takes
$\rho:\gamma\delta\la\alpha\beta$ and $\sigma:\delta B\la\Gamma\alpha$,
respectively, to
\[\xymatrix@R-10pt@C+20pt{
\delta B\ar[r]^-{\eta(\gamma\dashv\Gamma)} &\Gamma\gamma\delta B \ar[r]^-{\rho} & 
\Gamma\alpha\beta B\ar[r]^-{\e(\beta\dashv B)} & \Gamma\alpha\\
\gamma\delta\ar[r]^-{\eta(\beta\dashv B)} & \gamma\delta B\beta \ar[r]^-{\sigma} &
\gamma\Gamma\alpha\beta\ar[r]^-{\e(\gamma\dashv\Gamma)} & \alpha\beta
}\]
The standard terminology in category theory is that $\rho$ and $\s$ are
\emph{mates} of each other.
More explicitly: $\rho$ is the \emph{left mate} of $\s$, and $\s$ is
the \emph{right mate} of $\rho$.

As the reader may have guessed,
this terminology was invented by Australian category theorists.
\ermd

\con{C30.3}
Any commutative square of schemes
\[\xymatrix{
W \ar[r]^-u \ar[d]_f & X\ar[d]^g \\
Y\ar[r]^-v & Z  
}\]
gives rise to squares
\[\xymatrix@R-22pt@C+20pt{
\Dqc(W)&\ar[l]_-{\LL u^*} \Dqc(X) & & \Dqc(W)\ar[dddd]_-{\R f_*}&\ar[l]_-{\LL u^*} \Dqc(X) \ar[dddd]^-{\R g_*}\\
\quad  &\quad & &\quad &\quad  \\
\quad&\quad&\text{and} &\quad &\quad \\
\quad \ar@{=>}[ruu]^\tau &\quad& & \quad &\quad  \\
 \Dqc(Y)\ar[uuuu]^-{\LL f^*} &\ar[l]^-{\LL v^*}  \Dqc(Z)\ar[uuuu]_-{\LL g^*}& & \Dqc(Y) &
 \Dqc(Z)\ar[l]^-{\LL v^*} 
}\]
where $\tau$ is the canonical
isomorphism $\tau:\LL f^*\LL v^*\la \LL u^*\LL g^*$.
Reminder~\ref{R30.30303} yields a right mate for $\tau$, 
producing a natural transformation we will call
$\beta:\LL v^*\R g_*\la\R f_*\LL u^*$. The classical flat base-change
theorem
tells us

\sthm{SS30.3.1} The map $\beta$ is an isomorphism if the square is
cartesian---meaning a pullback square in the category of schemes---and 
if $v$ is flat.
\esthm
\econ

\con{C30.5}
As in Construction~\ref{C30.3}
assume given a commutative square of schemes
\[\xymatrix{
W \ar[r]^-u \ar@{}[dr]|{\diamondsuit}\ar[d]_f & X\ar[d]^g \\
Y\ar[r]^-v & Z  
}\]
If the base-change map  
$\R f_*\LL u^*\stackrel\beta\longleftarrow\LL v^*\R g_*$
is an isomophism, for example if we are in the situation of \ref{SS30.3.1},
we may apply Reminder~\ref{R30.30303} to the squares
\[\xymatrix@R-22pt@C+20pt{
\Dqc(W)\ar[dddd]_-{\R f_*}&\ar[l]_-{\LL u^*} \Dqc(X) \ar[dddd]^-{\R g_*}
 & & \Dqc(W)&\ar[l]_-{\LL u^*} \Dqc(X) \\
\quad \ar@{=>}[rdd]^{\beta^{-1}}  &\quad & &\quad &\quad  \\
\quad&\quad&\text{and} &\quad &\quad \\
\quad &\quad& & \quad &\quad  \\
\Dqc(Y) &
\Dqc(Z)\ar[l]^-{\LL v^*}
& &  \Dqc(Y)\ar[uuuu]^-{f^\times} &\ar[l]^-{\LL v^*}  \Dqc(Z)\ar[uuuu]_-{g^\times}
}\]
and produce a right mate for $\beta^{-1}$. We will
write this base-change map as 
$\Phi=\Phi(\diamondsuit):\LL u^*g^\times\la f^\times \LL v^*$.
And the relevant theorem tells us
\sthm{SS30.5.1} The map $\Phi(\diamondsuit)$ is an isomorphism if, in addition to the
hypotheses of~\ref{SS30.3.1}
already made to define $\Phi(\diamondsuit)$,
we assume that $g$ is proper and one of the following holds:
\be
\item
$f$ is of finite Tor-dimension.
\item
We restrict the functors to $\Dqcpl(Z)\subset\Dqc(Z)$; that is we only evaluate
the map on bounded-below complexes.
\ee
\esthm
\econ

\rmk{R30.7}
The assertion \ref{SS30.3.1} goes all the way back to Grothendieck in the
1950s; it's an easy consequence of the fact that one can compute $\R f_*$
using \u{Cech} complexes.

There is a (complicated) proof of \ref{SS30.5.1}(ii) 
in Hartshorne~\cite{Hartshorne66}, towards the end of the book.
The first short and formal proof of \ref{SS30.5.1}(ii) is due to
Verdier~\cite{Verdier68}. Over the
years there have been technical improvements, and the
strong version stated
in \ref{SS30.5.1}(i) is recent. The reader is referred to
\cite[Lemma~5.20]{Neeman13} for the proof.

For the applications in this article we do not need such refined forms
of \ref{SS30.5.1}. Verdier's old theorem 
 suffices.
\ermk

It's time to move on to the remaining business: the proofs of
Theorems~\ref{T0.13} and \ref{T0.19}.
Both theorems are assertions involving a certain natural map $\theta$;
the first step is to construct this $\theta$.
We're about to do this to show that it can be done formally,
using nothing more than the base-change map of Construction~\ref{C30.5},
the counit of some adjunction, and the Hochschild-Kostant-Rosenberg
Theorem. The reader willing to skip the construction should proceed
directly to
Remark~\ref{R25.9}.

\con{C30.9}
Suppose $f:X\la Y$ is flat and proper. Consider the diagram
\[\xymatrix@C+10pt@R+10pt{
  X\ar[r]^-\delta 
& X\times_Y^{} X\ar@{}[dr]|{\diamondsuit} \ar[r]^-{\pi_1^{}} \ar[d]_{\pi_2^{}} & X\ar[d]^f \\
& X\ar[r]^-f & Y  
}\]
where $\delta:X\la X\times_Y^{}X$ is the diagonal map.
Since $f$ is both flat and proper the hypotheses of \ref{SS30.5.1}
are satisfied, and
the map $\Phi(\diamondsuit)$ is an isomorphism.
We have the following composites, with the first
two defining the maps $\alpha$ and $\gamma$ that go into producing the third
\[\xymatrix@C+4pt@R-15pt{
f^\times\ar[r]_-\sim\ar@/^2pc/[rrr]^-{\alpha} &\LL\delta^* \LL\pi_1^* f^\times\ar[rr]^-{\LL\delta^*\Phi(\diamondsuit)}
&&\LL \delta^*\pi_2^\times \LL f^* \\
\ar@/_2pc/[rrr]^-{\gamma}\R\delta_*\ar[r]^-\sim &\R\delta_*\delta^\times\pi_2^\times\ar[rr]^-{\widetilde\e\pi_2^\times}
&&\pi_2^\times \\
& & &\\
\LL\delta^*\R\delta_*\LL f^*\ar[rr]^-{\LL\delta^*\gamma \LL f^*} &&
\LL\delta^*\pi_2^\times \LL f^*\ar[r]^-{\alpha^{-1}} & f^\times
}\]
The natural isomorphisms in
the first two composites are because $\pi_1^{}\delta=\pi_2^{}\delta=\id$,
hence $\LL\delta^*\LL\pi_1^*\cong\LL(\id)^*=\id=\id^\times\cong\delta^\times\pi_2^\times$.
Since the map $\Phi(\diamondsuit)$ is an isomorphism the first row
composes to an isomorphism, allowing us to form the third row.
Evaluating  the third row above at the object
$\co_Y^{}\in\Dqc(Y)$ we obtain the second and third maps
in the composite below defining $\zeta$
\[\xymatrix@C+20pt{
\LL\delta^*\R\delta_*\co_X^{} \ar@{=}[r] \ar@/^2pc/[rrr]^-\zeta&   
\LL\delta^*\R\delta_*\LL f^*\co_Y^{}\ar[r]_-{\LL\delta^*\gamma \LL f^*} &
\LL\delta^*\pi_2^\times \LL f^*\co_Y^{}\ar[r]_-{\alpha^{-1}} & f^\times\co_Y^{}
}\]
The object 
$\LL\delta^*\R\delta_*\co_X^{}$ is something we know and love---it is the
derived category version of the tensor product of $\co_X^{}$ with itself over
$\co_{X\times_Y^{}X}^{}$. Formal nonsense tells us that
$\LL\delta^*\R\delta_*\co_X^{}$ is a commutative monoid in the monoidal
category $\Dqc(X)$, hence its sheaf cohomology
$\ch^*(\LL\delta^*\R\delta_*\co_X^{})$ is a graded commutative ring.
There is an obvious ring homomorphism
\[\xymatrix@C+20pt{
\wedge^*\ch^{-1}[\LL\delta^*\R\delta_*\co_X^{}]\ar@{=}[r]&\wedge^*\Omega^1_{f}
\ar[r] & \ch^*(\LL\delta^*\R\delta_*\co_X^{})\ ,
}\]
from the exterior
algebra on $\ch^{-1}[\LL\delta^*\R\delta_*\co_X^{}]=\Omega^1_{f}$ to the ring
$\ch^*(\LL\delta^*\R\delta_*\co_X^{})$. So far we have only assumed $f$ flat
and proper.

Now assume $f$ is smooth and proper, of relative dimension $n$. The
Hochschild-Kostant-Rosenberg Theorem~\cite{Hochschild-Kostant-Rosenberg62}
(see also Lipman~\cite[Proposition~4.6.3]{Lipman87}) 
tells us that the homomorphism
of graded rings of the paragraph above
is an isomorphism.
In particular we deduce
\be
\item
The cohomology sheaves $\ch^i(\LL\delta^*\R\delta_*\co_X^{})$ vanish
for $i<-n$.
\item
We have constructed a natural map
\[\xymatrix@C+20pt{
\Omega^n_{f}\ar[r] &\wedge^n\ch^{-1}[\LL\delta^*\R\delta_*\co_X^{}]
\ar[r] & \ch^{-n}[\LL\delta^*\R\delta_*\co_X^{}]
}\]
\ee
Formal nonsense, about {\it t}--structures in triangulated categories,
tells us that the map of
cohomology sheaves
$\Omega^n_{f}\la\ch^{-n}[\LL\delta^*\R\delta_*\co_X^{}]$
can be realized as $\ch^{-n}$ of a unique morphism
$\psi:\Omega^n_{f}[n]\la\LL\delta^*\R\delta_*\co_X^{}$ in
the derived category $\Dqc(X)$. And the map
$\theta$ of Theorem~\ref{T0.13} is defined to be the composite
\[\xymatrix@C+20pt{
\Omega^n_{f}[n] \ar[r]^-\psi&   
\LL\delta^*\R\delta_*\co_X^{}\ar[r]^-{\zeta} & f^\times\co_Y^{}
}\]

The precise
version of Theorem~\ref{T0.13}  now says:
\econ

\thm{T25.1}
If $f$ is smooth and proper then 
the map $\theta$ of Construction~\ref{C30.9} is an isomorphism.
\ethm

\rmk{R25.7}
In passing we mentioned that the Verdier
version of \ref{SS30.5.1} is sufficient
for this paper---the reason is that
in the proof we will only evaluate the base-change maps
$\Phi(\diamondsuit)$
on the objects like $\co_Y^{}$ or $\pi^\times\co_X^{}$, which are bounded below.
These are the only objects that come up in the definition of the
map $\theta$ of Construction~\ref{C30.9}.

Lipman lectured about the
approach to the map $\theta$, presented in
Construction~\ref{C30.9}, already in the 1980s. But it only appeared
in print relatively recently, it may
be found in Alonso, Jerem\'{i}as and
Lipman~\cite[Example~2.4 and Proposition~2.4.2]{Alonso-Jeremias-Lipman12}.
The map $\theta$ has older avatars, for example in
Verdier~\cite{Verdier68}---although the fact that Lipman's and Verdier's maps
coincide was proved only recently, see \cite{Lipman-Neeman18}.
\ermk

\rmk{R25.9}
So far we have,
up
to one application of the Hochschild-Kostant-Rosenberg 
theorem,
 set up all the players entirely in formal nonsense fashion.
And the historical asides along the way tell the reader
that this much was known by the mid-1990s. Until now we haven't met anything
younger than two decades.

It remains to discuss the proofs of Theorems~\ref{T0.13} and \ref{T0.19},
and this is where there has been major progress in the last few years. We
open a new section for this.
\ermk

\subsection{The simple proofs discovered recently}
\label{SS30.5}

Now that we have defined the maps
occurring in Theorems~\ref{T0.13} and \ref{T0.19},
it remains to prove the claims
of the theorems---we need to show that the maps do
as the theorems
assert they do. This should be a local problem, but for
the longest time no one understood how to do this local computation
simply and elegantly. Recall the first paragraph of the Introduction: there
are two paths
to the foundations of the subject, and in this article we've been
following the one pioneered by Deligne and Verdier. The objection
to this approach has long been that it leads to a
theory where you can't compute anything. We've now reached the stage
where a computation is in order---it should hardly come as
a surprise that, until
very recently, this was the point where the trail we have been on
seemed to peter out, with the general direction appearing
blocked and impenetrable. 

We should probably explain, and for the purpose
of clarity let us narrow our attention to just Theorem~\ref{T0.13}.
Theorem~\ref{T0.19} is similar but a touch
more technical.

\rmk{R30.99995}
Since we will say almost nothing about the proofs of Theorem~\ref{T0.19},
old or new, we should in passing acknowledge its long history
and give references. The theorem was first
sketched in Hartshorne~\cite[pp.~398--400]{Hartshorne66}.
The case of varieties over a perfect field was completely worked out in
Lipman~\cite{Lipman84}. Lipman's main results were generalized in
H{\"u}bl and Sastry~\cite{Hubl-Sastry93}; see their
Residue Theorem (iii) on p.~752 and its generalization (iii) on p.~785.

The results cited in the paragraph above all came with complicated proofs.
The existence of a simple proof is a recent surprise, none of the experts---the
dozen of us---expected such a thing, it is part of the
exciting developments of the last few years. The reader can find
the simple proof in~\cite[Section~2]{Neeman13A}, and 
we will say a tiny bit more about this proof in \S\ref{SS35.1}.
\ermk

\str{ST30.380} Back to Theorem~\ref{T0.13}:
in Construction~\ref{C30.9} we defined a map
$\theta:\Omega^n_{f}[n]\la f^\times\co_Y^{}$, and
Theorem~\ref{T0.13} asserts that $\theta$ is an isomorphism.
Surely we should be able to prove this locally---a morphism in $\Dqc(X)$
is an isomorphism if and only if it restricts to an isomorphism on
an open affine cover.
Let us follow our noses and proceed the way such arguments usually work,
to see where the obstacle lies, and to pinpoint the insight
which uncovered an unconvoluted
path circumventing this hurdle.
\estr

\lem{L30.383}
Suppose we are given a commutative square of schemes
\[\xymatrix@R+5pt@C+20pt{
U\ar[d]_-{g}\ar[r]^-{u}\ar@{}[rd]|{\clubsuit} & X\ar[d]^-{f} \\
 V \ar[r]^-{v} & Y
}\]
where the maps $u$ and $v$ are open immersions.
Then there is a canonical natural isomorphism
$\LL u^*\Omega_{f}^n[n]\cong\Omega_{g}^n[n]$.
\elem

\prf
Obvious.
\eprf

\rmd{R30.385}
We are given
a morphism $f:X\la Y$, which we  assume smooth and proper.
We wish to show that the map
\[
\xymatrix@C+30pt{
\Omega_{f}^n[n]\ar[r]^-{\theta} & f^\times\co_Y^{}
}\]
of Construction~\ref{C30.9} is an isomorphism,
and the plan is to do this by studying it locally.
\ermd

\rdn{R30.28}
Let us begin by showing that the problem is local in $Y$. 
Suppose
that $v:V\la Y$ is an open immersion with $V$ affine. Next we
\be
\item
Form the pullback square of schemes
\[\xymatrix@R+5pt@C+20pt{
U\ar[d]_-{g}\ar[r]^-{u}\ar@{}[rd]|{\clubsuit} & X\ar[d]^-{f} \\
 V \ar[r]^-{v} & Y
}\]
It clearly suffices to show that, for every $v:V\la Y$
as above, the map
$\LL u^*\theta$ is an isomorphism. Our reduction is
about simplifying $\LL u^*\theta$.
\setcounter{enumiv}{\value{enumi}}
\ee
Let us pass from the square $(\clubsuit)$
to the square of derived categories
\[\xymatrix@R+5pt@C+20pt{
 \Dqc(U)&\ar[l]_-{\LL u^*} \Dqc(X)\\
 \Dqc(V)\ar[u]^-{g^\times} &\ar[l]_-{\LL v^*}  \Dqc(Y)\ar[u]_-{f^\times} 
}\]
The base-change map
$\Phi(\clubsuit):\LL u^*f^\times\la g^\times\LL v^*$ is
an isomorphism by \ref{SS30.5.1}. Applying to the object $\co_Y^{}\in\Dqc(Y)$,
and recalling that $\LL v^*\co_Y^{}=\co_V^{}$, we deduce
\be
\setcounter{enumi}{\value{enumiv}}
\item  We have produced an isomorphism
\[
\xymatrix{
\LL u^*f^\times\co_Y^{}\ar[r] &
g^\times\co_V^{}
}\]
\setcounter{enumiv}{\value{enumi}}
\ee
Next we apply Lemma~\ref{L30.383} to the square $(\clubsuit)$ in (i),
obtaining
\be
\setcounter{enumi}{\value{enumiv}}
\item  We have an isomorphism
\[\xymatrix@C+40pt{
\LL u^*\Omega_{f}^n[n]\ar[r]&  \Omega_{g}^n[n]
}\]
\setcounter{enumiv}{\value{enumi}}
\ee
We are trying to show that the map
\[\xymatrix@C+20pt{
\Omega_{f}^n[n]
\ar[r]^-{\theta} & f^\times\co_Y^{}
}\]
is an isomorphism, 
and we have already agreed in (i) that it suffices to check that
$\LL u^*\theta$ is an isomorphism.
By (ii) and (iii) above $\LL u^*\theta$ rewrites
as some map
\[\xymatrix@C+20pt{
\Omega_{g}^n[n]
\ar[r] & g^\times\co_V^{}
}\]
and it is an exercise in the definitions to check that this map is nothing
other than the $\theta$ corresponding to $g:U\la V$. In other words
we are reduced to proving Theorem~\ref{T0.13} in the case where $Y$
is affine.
\erdn

\rdn{R30.31}
Reduction~\ref{R30.28} tells us that it suffices to prove Theorem~\ref{T0.13}
for smooth, proper morphisms $f:X\la Y$ with $Y$ affine. If we follow the
usual yoga, the next step should be to reduce to the case where
$X$ is also affine.
We want to prove that the map $\theta$ of Construction~\ref{C30.9}
an isomorphism, and it certainly suffices
to show that $\LL u^*\theta$ is an isomorphism
for every open immersion $u:U\la X$ with $U$ affine. Choose therefore
an open immersion $u:U\la X$, and we would like to express $\LL u^*\theta$
in some form that renders it easily computable.
Now $\theta$ is a map
\[\xymatrix@C+20pt{
\Omega_{f}^n[n]
\ar[r] & f^\times\co_Y^{}
}\]
and we are embarking on a study of $\LL u^*\theta$, which is
a morphism
\[\xymatrix@C+20pt{
\LL u^*\Omega_{f}^n[n]
\ar[r]  &
\LL u^*f^\times\co_Y^{}
}\]
Let us begin by simplifying
$\LL u^*\Omega_{f}^n[n]$.

To this end we study the commutative square of schemes
\[\xymatrix@R+5pt@C+20pt{
U\ar[d]_-{fu}\ar[r]^-{u}\ar@{}[rd]|{\clubsuit} & X\ar[d]^-{f} \\
 Y \ar[r]^-{\id} & Y
}\]
Applying Lemma~\ref{L30.383}
produces an isomorphism 
\[\xymatrix@C+20pt{
\LL u^*\Omega_{f}^n[n]\ar[r]
&\Omega_{fu}^n[n]
}\]
In the usual jargon, our reduction so far tells us that the ``object
$\Omega_{f}^n[n]$ is local in $X$''.
\erdn

\cau{C30.35}
The simple-minded way to proceed doesn't work,
unfortunately 
the object $f^\times\co_Y^{}$ isn't local in $X$.
If we look at the composite $U\stackrel u\la X\stackrel f\la Y$,
then $\LL u^*f^\times\co_Y^{}$ is not in general
isomorphic to $(fu)^\times\co_Y^{}$. See Sketch~\ref{SK35.7} for a
counterexample.
\ecau

\rdn{R30.37}
Caution~\ref{C30.35} tells us what doesn't work when
we try to simplify $\LL u^*f^\times$, and it's now time
to present the way around the difficulty. To this end consider the following
commutative diagram of schemes
\[\xymatrix@C+10pt@R+5pt{
U \ar[r]^-u\ar[d]_\delta & X\ar@/^2pc/[ddrr]^-\id && \\
U\times_Y^{}U \ar@{}[dr]|{\clubsuit}\ar[r]^-\id\ar[d]_-\id & U\times_Y^{}U\ar[d]^-{u\times \id} &  & \\
U\times_Y^{}U \ar[r]^-{u\times\id} &
X\times_Y^{}U \ar@{}[dr]|{\diamondsuit} \ar[r]^-{\id\times u}\ar[d]_-{\pi_2^{}}
& X\times_Y^{}X \ar@{}[dr]|{\heartsuit} \ar[r]^-{\pi_1^{}}
\ar[d]_-{\pi_2^{}} & X\ar[d]^f \\  
& U \ar[r]^u &X \ar[r]^f & Y
}\] 
from which we obtain a diagram of derived categories
\[\xymatrix@C+10pt@R+5pt{
\Dqc(U) &\ar[l]_-{\LL u^*} \Dqc(X) && \\
\Dqc(U\times_Y^{}U)\ar[u]^{\LL\delta^*} &\ar[l]_-{\LL(\id)^*} \Dqc(U\times_Y^{}U) &  & \\
\Dqc(U\times_Y^{}U)\ar[u]^-{\id=\id^\times} &\ar[l]_-{\LL(u\times\id)^*} \ar[u]_-{(u\times\id)^\times} 
\Dqc(X\times_Y^{}U)& \ar[l]_-{\LL(\id\times u)^*}  \Dqc(X\times_Y^{}X) & \ar[l]_-{\LL\pi^*_1}
  \Dqc(X)\ar@/_2pc/[lluu]_-{\LL(\id)^*=\id} \\  
& \Dqc(U)\ar[u]_-{\pi_2^{\times}} &\ar[l]_-{\LL u^*} \Dqc(X)\ar[u]_-{\pi_2^{\times}} &\ar[l]_{\LL f^*}\ar[u]_{f^\times}  \Dqc(Y)
}\] 
We would like to simplify the expression $\LL u^*f^\times$, and
in the diagram above we have a subdiagram
which obviously commutes up to canonical natural isomorphism 
\[\xymatrix@C+10pt@R+5pt{
\Dqc(U) &\ar[l]_-{\LL u^*} \Dqc(X) && \\
\Dqc(U\times_Y^{}U)\ar[u]^{\LL\delta^*} & &  & \\
\Dqc(U\times_Y^{}U)\ar[u]^-{\id=\id^\times} &\ar[l]_-{\LL(u\times\id)^*} 
\Dqc(X\times_Y^{}U)& \ar[l]_-{\LL(\id\times u)^*}  \Dqc(X\times_Y^{}X) & \ar[l]_-{\LL\pi^*_1}
  \Dqc(X)\ar@/_2pc/[lluu]_-{\LL(\id)^*=\id} \\  
& & &\ar[u]_{f^\times}  \Dqc(Y)
}\] 
The squares $(\diamondsuit)$ and $(\heartsuit)$ in the commutative diagram
of schemes satisfy the hypotheses of \ref{SS30.5.1} and, up to
the natural 
isomorphisms induced by $\Phi(\diamondsuit)$ and $\Phi(\heartsuit)$
composed with the canonical natural isomorphism of the diagram above,
the following subdiagram must also commute
\[\xymatrix@C+10pt@R+5pt{
\Dqc(U) &\ar[l]_-{\LL u^*} \Dqc(X) && \\
\Dqc(U\times_Y^{}U)\ar[u]^{\LL\delta^*} & &  & \\
\Dqc(U\times_Y^{}U)\ar[u]^-{\id=\id^\times} &\ar[l]_-{\LL(u\times\id)^*} 
\Dqc(X\times_Y^{}U)&  &
\Dqc(X)\ar@/_2pc/[lluu]_-{\LL(\id)^*=\id} \\  
& \Dqc(U)\ar[u]_-{\pi_2^{\times}} &\ar[l]_-{\LL u^*} \Dqc(X)
&\ar[l]_{\LL f^*}\ar[u]_{f^\times}  \Dqc(Y)
}\]
Up until now everything is entirely classical.

Of course there is nothing to stop us from looking at the base-change map
of the square $(\clubsuit)$ in the large diagram
of schemes at the beginning
of this Reduction.
The square is cartesian, the
horizontal map $(u\times\id)$ is flat, and Construction~\ref{C30.5}
provides us with a morphism
$\Phi(\clubsuit):\LL(\id)^*(u\times\id)^\times\la \id^\times\LL(u\times\id)^*$,
or more simply
$\Phi(\clubsuit):(u\times\id)^\times\la \LL(u\times\id)^*$.
It might seem idiotic\footnote{In September 2016 the author presented
  the results (in a more technical incarnation)  at a seminar at
  the University of Utah. The audience consisted of algebraic
  geometers and commutative algebraists. And this expert audience burst out
  laughing when told that the study of the base-change
  map $\Phi(\clubsuit)$ is what underpins the recent progress---to the 
experts this was hilarious.}
to study $\Phi(\clubsuit)$,
after all the hypotheses of \ref{SS30.5.1}
don't hold, the vertical map on the right decidedly isn't proper.
And just in case the reader was wondering:
it's not just that the hypotheses of \ref{SS30.5.1}
aren't satisfied---neither is the conclusion, the
map $\Phi(\clubsuit)$ is known not to be an isomorphism in general.

The recent insight says
\sthm{ST30.39.1}
Consider the following extract from our large diagram
of derived categories
\[\xymatrix@C+10pt@R+5pt{
\Dqc(U) & \\
\Dqc(U\times_Y^{}U)\ar[u]^{\LL\delta^*} &\ar[l]_-{\LL(\id)^*} \Dqc(U\times_Y^{}U)\\
 \Dqc(U\times_Y^{}U)\ar[u]^-{\id^\times} &\ar[l]_-{\LL(u\times\id)^*} \ar[u]_-{(u\times\id)^\times} 
\Dqc(X\times_Y^{}U)
}\] 
Then the composites from bottom right to top left agree,
more precisely the map $\LL\delta^*\Phi(\clubsuit)$ gives
an isomorphism $\LL\delta^*(u\times\id)^\times\la \LL\delta^*\LL(u\times\id)^*$.
\esthm
With the aid of the isomorphism of \ref{ST30.39.1}
we deduce that, up to all the isomorphisms above,
the diagram
\[\xymatrix@C+10pt@R+5pt{
\Dqc(U) &\ar[l]_-{\LL u^*} \Dqc(X) && \\
\Dqc(U\times_Y^{}U)\ar[u]^{\LL\delta^*} &\ar[l]_-{\LL(\id)^*} \Dqc(U\times_Y^{}U) &  & \\
& \Dqc(X\times_Y^{}U)\ar[u]_{(u\times\id)^\times}&  &
\Dqc(X)\ar@/_2pc/[lluu]_-{\LL(\id)^*=\id} \\  
& \Dqc(U)\ar[u]_-{\pi_2^{\times}} &\ar[l]_-{\LL u^*} \Dqc(X)
&\ar[l]_{\LL f^*}\ar[u]_{f^\times}  \Dqc(Y)
}\]
also commutes. This simplifies to
\[\xymatrix@C-5pt@R+10pt{
& \Dqc(U) && \ar[ll]_-{\LL u^*} \Dqc(X) & \\
\Dqc(U\times_Y^{}U)\ar[ur]^{\LL\delta^*}   
& &\Dqc(U)\ar[ll]^-{\pi_2^{\times}} &&\ar[ll]^-{\LL (fu)^*} \ar[ul]_{f^\times}  \Dqc(Y)
}\]
and the punchline is that we have found an isomorphism of
the functor $\LL u^*f^\times$ with the
composite $\LL\delta^*\pi_2^{\times}\LL (fu)^*$,
and in that composite \emph{all the schemes are affine.}
\erdn

\rmk{R30.39}
Reductions~\ref{R30.31} and \ref{R30.37} transform the problem into one in
which all the schemes are affine. It remains to show that
\be
\item
  The affine problem is tractable, in other words we can do
  the computation needed in the affine case. The reader can find
  in~\cite[Theorem~4.2.4]{Iyengar-Lipman-Neeman13} a
  working out of what needs to be computed, and in
  \cite[Section~1]{Neeman13A} the computation is actually carried
  out. Below we present a rough outline, in
  Reminder~\ref{R30.40} and Computation~\ref{C30.41}.
\item
  We should also say something about the proof of \ref{ST30.39.1}, after
  all this is the crux of what's new
  to the approach of the current document. The expert is referred
  to~\cite[(2.3.5.1)]{Iyengar-Lipman-Neeman13} for a general
  result, which implies \ref{ST30.39.1} as a special case.
  In the interest of making the subject accessible to the
  non-expert we give a fairly complete, self-contained
  treatment in Sketch~\ref{SK30.45}
  below.
\ee
Let us first recall  
\ermk

\rmd{R30.40}
If $R$ is a ring let $\D(R)$ denote the (unbounded) derived category of cochain
complexes of $R$--modules. From \cite[Theorem~5.1]{Bokstedt-Neeman93} we
learn that, for any ring $R$, the canonical
functor $\D(R)\la\Dqc\big[\spec R\big]$
is an equivalence. Given a ring homomorphism $f:R\la S$ we have an induced
map of schemes, and by abuse of notation we will write it
$f:\spec S\la\spec R$. The diagram
\[\xymatrix@R+15pt{
\Dqc\big[\spec S\big]\ar[d]|-{\R f_*} \\
\Dqc\big[\spec R\big]\ar@<3ex>[u]^-{\LL f^*} \ar@<-3ex>[u]_-{f^\times}
}\]
of Reminder~\ref{R-1.-1} identifies with
\[\xymatrix@R+15pt{
\D(S)\ar[d]|-{f_*} \\
\D(R)\ar@<3ex>[u]^-{\LL f^*} \ar@<-3ex>[u]_-{f^\times}
}\]
where $f_*:\D(S)\la\D(R)$ is the forgetful
functor---it takes an object of $\D(S)$, that
is a cochain complex of $S$--modules, to itself viewed as a complex of
$R$--modules. As $f_*$ is exact there is no need to derive it,
we have $\R f_*=f_*$. With the notation as in Example~\ref{E24.7}
the functor $\LL f^*$, being the left adjoint of $f_*$, is
given by the formula $\LL f^*(-)=S\oo_R^\LL(-)$,
while $f^\times$, being the right adjoint of $f_*$, is the functor
$f^\times(-)=\RHom_R^{}(S,-)$. And the units and counits of the
adjunctions are all explicit.
\ermd

\cmp{C30.41}
In view of Reminder~\ref{R30.40}, achieving Remark~\ref{R30.39}(i)
has to be straightforward---it's just a matter of untangling the definitions
and then doing a computation. 
We are given affine schemes
$U$ and $Y$, hence we may write $Y=\spec R$ and $U=\spec S$. The morphism
of schemes $fu:U\la Y$ corresponds to a ring homomorphism $\s:R\la S$,
and we have an isomorphism $ U\times_Y^{}U=\text{Spec}(\se)$ with
$\se=S\oo_R^{}S$. Reductions~\ref{R30.31} and \ref{R30.37}
produce isomophisms in the category $\Dqc(U)$
\begin{eqnarray*}
 \LL u^*\Omega_{f}^n[n]&\cong&
  \Omega_{fu}^n[n] \\
  \LL u^*f^\times\co_Y^{}&\cong&
  \LL\delta^*\pi_2^{\times}\LL (fu)^*\co_Y^{}
\end{eqnarray*}
Using the descriptions on the right-hand-side,
one checks that
the equivalence $\Dqc(U)\cong\D(S)$ takes
these objects, respectively,
to $\Tor^\se_n(S,S)[n]$ and $S\oo_\se^\LL\RHom_R^{}(S,S)$.
And the map $\LL u^*\theta$
is nothing other than the composite
\[\xymatrix@C+20pt{
\Tor^\se_n(S,S)[n]\ar[r] &S\oo_\se^\LL S\ar[r]^-{\id\oo_{\se}^\LL I} &S\oo_\se^\LL\text{RHom}_R^{}(S,S)
}\]
where $I:S\la\text{RHom}_R^{}(S,S)$ is the obvious inclusion. It remains
to show that, for $S$ smooth over $R$ of relative dimension $n$, the
composite above is an isomorphism. The reader can find the computation
in \cite[Section~1]{Neeman13A}.
\ecmp

\skt{SK30.45}
It remains to deliver on the promise of Remark~\ref{R30.39}(ii),
we should say something about the proof of \ref{ST30.39.1}.
The argument below is reasonably detailed---it
may safely be skipped, the reader
should feel free to proceed directly to Section~\ref{S35}.

We remind the reader: we consider the diagram
\[\xymatrix@C+10pt@R+5pt{
\Dqc(U) & 
\Dqc(U\times_Y^{}U)\ar[l]_-{\LL\delta^*} &\ar[l]_-{\LL(\id)^*} \Dqc(U\times_Y^{}U)\\
& \Dqc(U\times_Y^{}U)\ar[u]^-{\id^\times} &\ar[l]_-{\LL(u\times\id)^*} \ar[u]_-{(u\times\id)^\times} 
\Dqc(X\times_Y^{}U)
}\]
and the assertion is that the functor
$\LL\delta^*$ takes the base-change map
$\Phi(\clubsuit):(u\times\id)^\times\la \LL(u\times\id)^*$ to an
isomorphism. To simplify the notation we will write $v$ for
the map $(u\times\id):U\times_Y^{}U\la X\times_Y^{}U$.

Now the diagonal map
$\delta:U\la U\times_Y^{}U$ is a closed immersion, hence the map
$\R\delta_*=\delta_*$ is conservative---it suffices to prove
that $\R\delta_*\LL\delta^*$ takes $\Phi(\clubsuit)$ to an isomorphism.
But we have isomorphisms
\[
\R\delta_*\LL\delta^*(-)\quad\cong\quad
\R\delta_*\big[\LL\delta^*(-)\oo^\LL\co_U^{}\big]\quad\cong\quad
(-)\oo^\LL\R\delta_*\co_U^{}
\]
where the second is the isomorphism $p^{-1}$ of the projection formula,
see Sketch~\ref{SK30.3}.
Consider therefore the two full subcategories
of $\Dqc(U\times_Y^{}U)$
given by
\[
  \D_{\mathbf{qc},\Delta}^{}(U\times_Y^{}U)\eq\left\{\mathfrak{s}\in\Dqc(U\times_Y^{}U)\left|
  \begin{array}{c} \LL i^*\mathfrak{s}=0\text{ for all }\\
i:\spec K\la U\times_Y^{}U-\Delta\\
    \text{where $K$ is a field and where}\\
    \Delta\subset U\times_Y^{}U\text{ is the
      diagonal}\end{array}
\right.\right\}
\]
\[
\cs\eq\left\{\mathfrak{s}\in\Dqc(U\times_Y^{}U)\left|\begin{array}{c}
\text{the functor }(-)\oo^\LL \mathfrak{s}\\
\text{ takes }\Phi(\clubsuit)\text{ to an isomorphism}
\end{array}\right.\right\}
\]
The object $\R\delta_*\co_U^{}$ clearly belongs to
$\D_{\mathbf{qc},\Delta}^{}(U\times_Y^{}U)$, and we wish to show that it
belongs to $\cs$. It certainly suffices to prove
that $\D_{\mathbf{qc},\Delta}^{}(U\times_Y^{}U)$ is
contained in $\cs$.

But both $\cs$ and $\D_{\mathbf{qc},\Delta}^{}(U\times_Y^{}U)$ are
localizing tensor ideals, and
\cite[Corollary~3.4]{Neeman92B} tells us that, as a localizing tensor ideal,
$\D_{\mathbf{qc},\Delta}^{}(U\times_Y^{}U)$ is generated by the perfect
complexes inside it\footnote{Note that $U\times_Y^{}U$ is affine, so for
us the old version in~\cite{Neeman92B} suffices. We should mention that
\cite{Neeman92B} builds on earlier papers by Hopkins~\cite{Hopkins85}
and Thomason and Trobaugh~\cite{ThomTro}.
The reader can find
later improvements in: the union
of Thomason~\cite[Lemma~3.4]{Thomason97} and
Alonso, Jerem{\'\i}as and
Souto~\cite[Corollary~4.11 and Theorem~4.12]{Alonso-Jeremias-Souto04}
generalize the result to all noetherian schemes, 
while Balmer and Favi~\cite{Balmer-Favi11} give a
formal generalization to the world of
tensor triangulated categories.}. It suffices to show that the perfect
complexes in
$\D_{\mathbf{qc},\Delta}^{}(U\times_Y^{}U)$ all belong to $\cs$.
In other words: it suffices to prove that,
for every perfect complex $P$ supported on the
diagonal, the functor $(-)\oo^\LL P$ takes 
$\Phi(\clubsuit)$ to an isomorphism. 

The morphism $v:U\times_Y^{}U\la X\times_Y^{}U$ is an open immersion,
hence the counit of adjunction $\e:\LL v^*\R v_*\la\id$ is an isomorphism.
Consequently $\R v_*$ is fully faithful, and we have an isomorphism
$\LL v^*\R v_*P\cong P$. It therefore suffices to show that, for
every perfect complex $P\in\Dqc(U\times_Y^{}U)$ supported on the diagonal,
the functor
\[
\R v_*\big[(-)\oo^\LL P\big]\quad\cong\quad
\R v_*\big[(-)\oo^\LL\LL v^*\R v_*P\big]\quad\cong\quad
\R v_*(-)\oo^\LL \R v_*P
\]
takes $\Phi(\clubsuit)$ to an isomorphism, or to put it differently
the functor $(-)\oo^\LL\R v_*P$ takes $\R v_*\Phi(\clubsuit)$ to
an isomorphism. Now let $P\in\Dqc(U\times_Y^{}U)$
be a perfect complex
supported on the diagonal and let $\Gamma\subset X\times_Y^{}U$
be the graph 
of the map $u:U\la X$. The following is a cartesian square of open immersions
\[\xymatrix@C+40pt{
U\times_Y^{}U-\Delta \ar[r]^-\alpha \ar[d]_-\beta & U\times_Y^{}U\ar[d]^-v\\
X\times_Y^{}U-\Gamma \ar[r]^-\gamma & X\times_Y^{}U
}\]
By \ref{SS30.3.1} we have the first isomorphism below
\[
\LL\gamma^*\R v_* P\quad\cong\quad \R\beta_*\LL\alpha^* P\quad\cong\quad 0\ ,
\]
where the second isomorphism is because $P$ is supported on $\Delta$ and
hence $\LL\alpha^*P=0$.
Thus both $\LL\gamma^*\R v_* P=0$ and $\LL v^*\R v_* P\cong P$ are
perfect complexes, and as $U\times_Y^{}U$ and $X\times_Y^{}U-\Gamma$ form
an open cover for $X\times_Y^{}U$ we deduce
\be
\item
If $P\in\Dqc(U\times_Y^{}U)$ is a perfect complex supported on the  
diagonal then 
$\R v_*P$ is a perfect complex
on $X\times_Y^{}U$, supported on the graph $\Gamma\subset X\times_Y^{}U$
of the map $u:U\la X$.
\setcounter{enumiv}{\value{enumi}}
\ee
The reader can amuse herself by proving
\be
\setcounter{enumi}{\value{enumiv}}
\item
The map $\Phi(\clubsuit):v^\times\la \LL v^*$ is taken by $\R v_*$  to
the composite 
\[\xymatrix@C+20pt{
\R v_*v^\times\ar[r]^-{\wt\e} & \id \ar[r]^-{\eta} & \R v_*\LL v^*
}\]
where $\wt\e:\R v_*v^\times\la\id$ is the counit of the adjunction
$\R v_*\dashv v^\times$, while $\eta:\id\la\R v_*\LL v^*$ is the unit
of the adjunction $\LL v^*\dashv\R v_*$.
\setcounter{enumiv}{\value{enumi}}
\ee
Let $Q\in\Dqc(X\times_Y^{}U)$ be a perfect complex
supported on $\Gamma$.
Then its dual $Q^\vee$ is also a perfect
complex supported on $\Gamma$, and
we have $(-)\oo^\LL Q\cong\RHHom(Q^\vee,-)$.
From (i) and (ii) above it suffices to prove
\be
\setcounter{enumi}{\value{enumiv}}
\item
For all perfect
complexes $Q\in\Dqc(X\times_Y^{}U)$ supported on $\Gamma$,
the functor
$(-)\oo^\LL Q$ takes the map 
$\eta$ of (ii) to an isomorphism.
\item
For all perfect
complexes $Q\in\Dqc(X\times_Y^{}U)$ supported on $\Gamma$,
the functor
$\RHHom(Q,-)$
takes the map $\wt\e$ of (ii) to an isomorphism.
\setcounter{enumiv}{\value{enumi}}
\ee
To establish (iv) it suffices, by \cite[Lemma~3.2]{Neeman96},
to show that for all pairs of perfect complexes
$C,Q\in\Dqc(X\times_Y^{}U)$, with
$Q$ supported on $\Gamma$, the functor
$\Hom(C,-)$ takes $\RHHom(Q,\wt\e)$ to an isomorphism.
Now observe the isomorphism of functors
\[
\Hom\big[C\,\,,\,\,\RHHom(Q,-)\big]
\quad\cong\quad \Hom(C\oo^\LL Q\,\,,\,\,-)
\]
As $C$ and $Q$ are both perfect and $Q$ is supported on $\Gamma$,
the complex $C\oo^\LL Q$ is perfect and is supported on $\Gamma$.
Hence (iv) would follow from
\be
\setcounter{enumi}{\value{enumiv}}
\item
For all perfect complexes $Q\in\Dqc(X\times_Y^{}U)$
supported on $\Gamma$, the functor
$\Hom(Q,-)$
takes the map $\wt\e$ of (ii) to an isomorphism.
\setcounter{enumiv}{\value{enumi}}
\ee

Let $Q\in\Dqc(X\times_Y^{}U)$ be a perfect complex
supported on $\Gamma$. Then $\LL v^*Q$ is a perfect complex
supported on $\Delta$. The map $\eta:Q\la\R v_*\LL v^*Q$ is an
isomorphism on the open set $U\times_Y^{}U\subset X\times_Y^{}U$,
and is an isomorphism on
$(X\times_Y^{}U) -\Gamma$ because both $Q$ and  [by (i)] $\R v_*\LL v^*Q$
vanish outside $\Gamma$. Since the open
sets $U\times_Y^{}U$ and $(X\times_Y^{}U) -\Gamma$ cover
$X\times_Y^{}U$ it follows that $\eta:Q\la\R v_*\LL v^*Q$
is an isomorphism.
Putting $A=\LL v^*Q$ we deduce that
(iii) and (v) would follow from
\be
\setcounter{enumi}{\value{enumiv}}
\item
For all $A\in\Dqc(U\times_Y^{}U)$, the functor
$(-)\oo^\LL \R v_*A$ takes
$\eta$ to an isomorphism.
\item
For all $A\in\Dqc(U\times_Y^{}U)$, the functor
$\Hom(\R v_*A\,,\,-)$
takes $\wt\e$ to an isomorphism.
\setcounter{enumiv}{\value{enumi}}
\ee
To see (vi) observe the isomorphism of the projection
formula
\[
(-)\oo^\LL \R v_*A\quad\cong\quad\R v_*\big[\LL v^*(-)\oo^\LL A\big]
\]
Since the functor $\LL v^*$ takes $\eta$ to an ismorphism so does
the right-hand-side above, and hence also the left-hand-side.

To see
(vii) observe the commutative square
\[\xymatrix@C+40pt{
\Hom\big(\R v_*\LL v^*(-)\,\,,\,\,?\big)  
\ar[r]^-\sim  \ar[d]_-{\Hom(\eta,?)}&
\Hom\big(-\,\,,\,\,\R v_*v^\times(?)\big)  \ar[d]^-{\Hom(-,\wt\e)}\\
\Hom(-,?) \ar@{=}[r] & \Hom(-,?)
}\]
We wish to show that the vertical map on
the right is an isomorphism when evaluated at $(-)=\R v_*A$, and
the vertical map on the left makes it clear, after all
$\eta\R v_*:\R v_*\la\R v_*\LL v^*\R v_*$ is an isomorphism.
\eskt

\section{Why did it take so long?}
\label{S35}

In some sense the ingredients of the proof were available already in the
1960s, but back then no one thought of applying the tools of homotopy
theory---for example Brown representability---to problems in algebraic
geometry. The methods employed in the classical proofs are fundamentally
unsuited for the approach presented here. To mention just one facet: in
the argument given here we relied heavily on the full power of the
derived tensor product. The pre-1990 literature on Grothendieck duality
all worked in the bounded-below derived category, where the derived tensor
product exists only under strong restrictions and is next to useless.

That said, with the exception of \ref{ST30.39.1}
the ingredients of the argument were all available by the mid-1990s.
And in Sketch~\ref{SK30.45} the reader learned that the proof
of \ref{ST30.39.1} could also have been given two decades ago. So why did we
fail to see this?

\subsection{The strangeness of the argument}
\label{SS35.1} The reader should appreciate the bizarreness of looking
at the base-change map in \ref{ST30.39.1}. It might help to elaborate
a little by sketching a related and slightly easier computation.

\skt{SK35.7}
Let $k$ be a field, let
$f:X\la Y$ be the projection $f:\pp^n_k\la\spec k$,
and let $u:U\la X$ be the open immersion $\ak^n_k\la\pp^n_k$. In
Caution~\ref{C30.35} we warned the reader that the functor
$u^\times f^\times$ is very different from $\LL u^*f^\times$. Let
us work out just how different by evaluating on $\co_Y^{}$.
Theorem~\ref{T0.13} gives an isomorphism
$\theta:\Omega^n_{f}[n]\la f^\times\co_Y^{}$, which
means that
$\LL u^*f^\times\co_Y^{}\cong\LL u^*\Omega^n_{f}[n]\cong\Omega^n_{\ak_k^n}[n]\cong
\co_U^{}[n]$,
after all the canonical bundle of $U=\ak^n_k$ is trivial. Under the equivalence
$\D(S)\cong\Dqc(U)$ of Reminder~\ref{R30.40} and Computation~\ref{C30.41}
the object  $\LL u^*f^\times\co_Y^{}\in\Dqc(U)$ identifies with $S[n]\in\D(S)$,
that is the cochain complex obtained by placing the $S$--module $S$
in dimension $-n$.

Now let's compute $u^\times f^\times\co_Y^{}\cong(fu)^\times\co_Y^{}$.
The morphism $fu:U\la Y$ is a map of affine schemes, and in
Reminder~\ref{R30.40} we noted that the computation
of $(fu)^\times$ can be carried over to $(fu)^\times:\D(R)\la\D(S)$ and
is given explicitly by the formula $(fu)^\times(-)=\RHom_R^{}(S,-)$.
In our case $R$ is the field $k$, $S=k[x_1^{},\ldots,x_n^{}]$
is the polynomial ring, and $(fu)^\times\co_Y^{}$ computes
to be $\RHom_k(S,k)=\Hom_k(S,k)$, which is a gigantic injective
$S$--module placed in degree zero. The reader can check
\cite{Neeman13B} to see just how gargatuan this injective module is.
It turns out to depend on the cardinality of $k$.

Now consider the cartesian square
\[\xymatrix{
\ak^n_k \ar[r]^\id\ar[d]_\id & \ak^n_k\ar[d]^u \\
\ak^n_k\ar[r]^u & \pp^n_k
}\]
The bottom horizontal map is flat, so there is a base-change map
$\Phi:\LL(\id)^*u^\times\la\id^\times\LL u^*$; more
simply we can write it as $\Phi:u^\times\la\LL u^*$. Since $\Phi$ is
a natural tranformation between two functors
$\xymatrix{\Dqc(\pp^n_k)\ar@<0.3ex>[r]\ar@<-0.3ex>[r]&\Dqc(\ak^n_k) }$
we may evaluate it at the object 
$f^\times\co_{\spec k}^{}\in\Dqc(\pp^n_k)$, producing
a morphism
$\psi=\Phi_{f^\times\co_{\spec k}^{}}^{}:u^\times f^\times\co_{\spec k}^{}
\la\LL u^*f^\times\co_{\spec k}^{}$
in the category $\Dqc(\ak^n_k)$. The equivalence
$\Dqc(\ak^n_k)\cong\D(S)$ must take it to a morphism
\[\xymatrix@C+20pt{
\Hom_k(S,k)\ar[r]^-{\psi} & S[n]
}\]
Doesn't it seem absurd to study this map?

But this is exactly what we do in the simple, recent
proof of Theorem~\ref{T0.19},
which we haven't discussed in this article. The next paragraph gives
a quick sketch---the non-experts may wish to skip ahead to \S\ref{SS35.295}.

Let $W\subset X$ be as in the statement of Theorem~\ref{T0.19}---that
is $W$ is the set of closed points in $X=\pp^n_k$,
and hence $U\cap W\subset U$ the set of closed points in $U=\ak^n_k$.
The recent proof of Theorem~\ref{T0.19} hinges on the
observation that the functor $\Gamma_{U\cap W}^{}$ takes
the map $\psi$ above to an
isomorphism. In Example~\ref{E0.17} we told the
reader how to compute $\Gamma_W^{}f^\times\co_{Y}^{}$, and
the analogous recipe gives that
$\Gamma_{U\cap W}^{}\LL u^*f^\times\co_{Y}^{}
\cong\LL u^*\Gamma_W^{}f^\times\co_Y^{}$\footnote{The isomorphism comes from
Remark~\ref{R-1.999}}
can be computed by forming a minimal injective resolution for
$S[n]$, that is the complex
\[\xymatrix{
  0\ar[r] & I^{-n}\ar[r] &I^{-n+1}\ar[r] &
  \cdots\ar[r] & I^{-1}\ar[r] & I^0\ar[r] & 0
}\]
and putting $I^0=\Gamma_{U\cap W}^{}\LL u^*f^\times\co_{Y}^{}$.
The morphism $\psi$ in the category $\D(S)$ may be represented by a
cochain map of injective resolutions
\[\xymatrix{
  0\ar[r] & 0\ar[r]\ar[d] &0\ar[r]\ar[d] &
  \cdots\ar[r] & 0\ar[r]\ar[d] & \Hom_k(S,k)\ar[r]\ar[d]^{\wt\psi} & 0 \\
   0\ar[r] & I^{-n}\ar[r] &I^{-n+1}\ar[r] &
  \cdots\ar[r] & I^{-1}\ar[r] & I^0\ar[r] & 0
}\]
and the functor $\Gamma_{U\cap W}^{}$ takes this to a map
$\Gamma_{U\cap W}^{}(\wt\psi):\Gamma_{U\cap W}^{}\Hom_k(S,k)\la I^0$.
We haven't yet told the reader how to compute $\Gamma_{U\cap W}^{}\Hom_k(S,k)$;
the formula is that $\Gamma_{U\cap W}^{}\Hom_k(S,k)\subset\Hom_k(S,k)$
is the submodule of all elements supported at the closed
points. That is: an element $e\in\Hom_k(S,k)$
belongs to the submodule $\Gamma_{U\cap W}^{}\Hom_k(S,k)$ if there exist
a finite number of maximal ideals $\m_1^{},\m_2^{},\ldots,\m_\ell^{}\subset S$
and an integer $N>0$ so that the ideal $\m_1^{N}\m_2^{N}\cdots\m_\ell^N$
annihilates $e$. So the assertion is that the obvious composite
\[\xymatrix@C+20pt{
\Gamma_{U\cap W}^{}\Hom_k(S,k) \ar@{^{(}->}[r] &\Hom_k(S,k)\ar[r]^-\psi &  I^0
}\]
is an isomorphism of $S$--modules. In \cite[Section~2]{Neeman13A} the reader
can see how this is used in the proof of Theorem~\ref{T0.19}.
\eskt

\subsection{The historical block}
\label{SS35.295} The fact that the key new idea is so outlandish is only part
of our excuse for taking so long. There were also the historical circumstances.

The foundations of Grothendieck duality was a lively, active field from about
the mid-1960s until well into the 1970s. And then the interest
gradually faded.
The mathematical community accepted that the foundations were complicated,
and stopped thinking about it. As an editor of one
journal put it, when
rejecting a recent paper of mine on the subject:
``\ldots your paper will be read by only the hardcore people\ldots
people have moved on\ldots''.
It's pretty accurate to say that, in the last three decades,
there have been two groups of people who have studied the foundations
of Grothendieck duality, in the old-fashioned setting of
classical, ordinary schemes: Lipman
and his students and collaborators on one side, and Yekutieli and his students
and collaborators on the other. It's not the whole story, one can point
to exceptions such as
Conrad~\cite{Conrad00}, but it is a good approximation of the truth.
In this I must count as one of Lipman's collaborators---every few years
we run into each other, and he persuades me to return to the problem.

Lipman's approach has long been orthogonal to Yekutieli's: Yekutieli accepted
the Grothendieck formalism, while Lipman has largely
followed the approach of Deligne and Verdier.
In this survey we've said next to nothing about the Grothendieck
angle, making any comparison difficult---one could fairly say
it takes an expert to recognize that the two subjects are truly one
and the same. In fact: checking the details, that is verifying that the
maps defined in the two theories agree, is often nontrivial. Just ask
Lipman---he's probably the person who has tried hardest.

To put it in a nutshell: the story isn't that hundreds of people were
feverishly working away
at the problem, and this multitude overlooked the obvious for twenty-some years.
What actually transpired is that at most a dozen of us were still studying
foundational questions, half of us were exploring 
what turned out to be
the wrong continent, and the solution, when it ultimately came, scored pretty
high on the weirdness scale.

\subsection{What finally woke us up}
\label{SS35.5}

Back to the recent progress: the
stimulus which, at long last, nudged us into probing in the right direction
came from the work of Avramov and
Iyengar~\cite{Avramov-Iyengar08}, and later
Avramov, Iyengar, Lipman and Nayak~\cite{Avramov-Iyengar-Lipman-Nayak10}.
They discovered puzzling formulas, and since then there has been a
proliferation---the reader can find increasingly general formulas
in~\cite{Iyengar-Lipman-Neeman13,NeemanTIFR}.
Let us present one example.
Suppose we are given composable morphisms of schemes
$U\stackrel u\la X\stackrel f\la Y$, and assume that $Y=\spec R$ and $U=\spec S$
are affine. Suppose the map $u$ is an open immersion, the map $f$ is
proper and the composite $fu:U\la Y$ is flat\footnote{The flatness
  of $fu$ isn't
  crucial, it suffices
  for $fu$ to be of finite Tor-dimension---the
  formula still
holds, it is a special case of Avramov, Iyengar, Lipman and
Nayak~\cite[(4.6.1)]{Avramov-Iyengar-Lipman-Nayak10}. But in the
finite-Tor-dimension generality
we don't have
a proof that's elementary. We will return to this in Problem~\ref{P95.9}.}.
The formula we happened to choose says that, for any $N\in\D(R)$, 
\[
\LL u^*f^\times N\quad\cong\quad S\oo_\se^\LL\RHom_R^{}\big(S,S\oo_R^\LL N\big)\ ,
\]
where $\se=S\oo_R^{}S$. Where on earth did this come from?

The original presentation, in the
articles~\cite{Avramov-Iyengar08,Avramov-Iyengar-Lipman-Nayak10}, touted
the left-hand-side as a great simplification of the right-hand-side.
Even the name given to the formulas---the ``Reduction Formulas''---reflects
this
perspective.
The formulas were first proved using the full strength
of the existing theory of Grothendieck duality.
It was not until \cite{Iyengar-Lipman-Neeman13}
that we came up with an elementary proof of the formulas\footnote{
  As exposed in \cite{Iyengar-Lipman-Neeman13} the proof doesn't seem
  elementary---the article \cite{Iyengar-Lipman-Neeman13}
  was written for an expert audience. But the formula
  is a minor variant of Reduction~\ref{R30.37} and, as presented in 
  this document, the proof is manifestly
  elementary. And the truth is that, modulo peeling away
  the generality in \cite{Iyengar-Lipman-Neeman13} and dusting off superfluous
  fluff, the proof here is identical to the proof there.},
and not
until \cite{Neeman13A} that the right-hand-side became a computational
tool for working out what's on the left.

What can I say: we were slow to see the light.

\section{Future directions}
\label{S95}

For the last three decades Grothendieck duality has been a small, niche subject,
with
only a handful of dedicated practitioners.
Granted that, the reader might well wonder what this section could
possibly be about. What conceivable future can there be in a moribund, small
field, long abandoned by the hordes?

There are two customary causes for the waning of a subject: it may die
because the important questions have all been satisfactorily answered,
or else because it hits a brick wall, and no one has any idea how to
advance. In both cases a rebirth is possible, but it takes something of
an earthquake. There needs to be a startling new development, a major
new insight, or vital new questions that come up.

Grothendieck duality is an example of a subject that died because people
got stuck---there were plenty of simple,
natural questions left, but no really good ideas on how to
tackle them. When we view the recent progress against this background, it can
only be a matter of time before the field picks up again---if
nothing else, the developments offer a radically new slant on
the what's known. In this
section we will sketch some of the obvious problems that now seem
within reach. But before the open questions we need to brush up on
facts that are known but haven't been covered yet---until now
we have made a conscious effort to be minimal in our use of category
theory, if a category or a functor was dispensable we omitted it.

\subsection{Assorted background material}
\label{SS95.1}

We have told the reader how to compute the functor $f^\times$ when
$f$ is smooth and proper. There is another classical situation in
which the computation of $f^\times$ is understood, let us recall.

\rmk{R95.3.1}
The category $\Dqc(X)$ is monoidal, it has a tensor product. This
tensor product has a right adjoint: there is a functor
$\RHHom_{\Dqc(X)^{}}^{}(-,?)$ and an isomorphism, natural in everything
in sight,
\[
\Hom\big[A\oo^\LL_{}B\,,\,C\big]\quad\cong\quad
\Hom\big[A\,,\,\RHHom_{\Dqc(X)^{}}^{}(B,C)\big]
\]
One way to construct it is to fix $B$ and note that, since
the functor $(-)\oo^\LL_{}B$ respects coproducts, Brown representability
gives a
right adjoint $\RHHom_{\Dqc(X)^{}}^{}(B,-)$. Now let $f:X\la Y$ be a
morphism of schemes and 
suppose $A,C$ are objects of $\Dqc(Y)$. We remind the reader of the following
string of isomorphisms
\begin{eqnarray*}
\Hom\big[A\,,\,\R f_*f^\times C\big]
& \cong & \Hom\big[\LL f^*A \,,\,f^\times C\big]\\
& \cong & \Hom\big[\LL f^*A\oo_{}^\LL\co_X^{} \,,\,f^\times C\big]\\
  & \cong & \Hom\big[\R f_*(\LL f^*A\oo_{}^\LL\co_X^{} )\,,\,C\big]\\
  & \cong & \Hom\big[A\oo_{}^\LL \R f_*\co_X^{}\,,\,C\big]\\
  & \cong & \Hom\big[A\,,\,\RHHom_{\Dqc(Y)^{}}^{}(\R f_*\co_X^{},C)\big]
\end{eqnarray*}
where the second isomorphism is because $\co_X^{}$ is the unit of
the tensor, the fourth comes from the projection formula,
and the others are all by adjunction.
Yoneda tells us that we have produced an isomorphism, natural
in $C$,
\[
\R f_*f^\times C\quad\cong\quad
\RHHom_{\Dqc(Y)^{}}^{}(\R f_*\co_X^{},C)
\]
When $f$ is an affine morphism this looks great---after all for
an 
affine morphism $f$ the functor $\R f_*$ is informative,
especially when we view it as a functor
from $\Dqc(X)$ to $\Dqc(\Mod{f_*\co_X^{}})$.

The problem is that, in general, the expression
$\RHHom_{\Dqc(Y)^{}}^{}(\R f_*\co_X^{},C)$ isn't overly computable.
What is frequently far more amenable to calculation is 
$\RHHom_{\D(Y)^{}}^{}(\R f_*\co_X^{},C)$, which is related. 
That is: we compute the internal
Hom not in the category $\Dqc(Y)$ but in the larger category
$\D(Y)$, whose objects are all
cochain complexes of sheaves of $\co_Y^{}$--modules.
This functor has the property that the cohomology sheaves
$\ch^i\big[\RHHom_{\D(Y)^{}}^{}(\R f_*\co_X^{},C)\big]$
are what one would hope for.
This means: for every open immersion
$u:U\hookrightarrow X$ consider
the $\Gamma(U,\co_Y^{})$--module
$\Hom_{\Dqc(U)^{}}^{}\big(\LL u^*\R f_*\co_X^{},\LL u^*C[i]\big)$. This comes
with obvious restriction maps, elevating the construction to a
presheaf 
of $\co_Y^{}$--modules on
$Y$. 
And the sheaf $\ch^i\big[\RHHom_{\D(Y)^{}}^{}(\R f_*\co_X^{},C)\big]$
turns out to be the sheafification of this presheaf. 
Unfortunately
in
order to pass, from the computable $\RHHom_{\D(Y)^{}}^{}(\R f_*\co_X^{},C)$
to the desired 
$\RHHom_{\Dqc(Y)^{}}^{}(\R f_*\co_X^{},C)\big]\cong\R
f_*f^\times C$, 
one
has to ``quasicoherate''---not the world's most transparent
process.\footnote{In Appendix~\ref{A97} we compute a simple example
  to illustrate the point.}    

There are cases where the computable Hom already has quasicoherent
cohomology, in which case the two functors agree and the
quasicoherator does nothing. Since we're making the assumption that
$f$ is an affine map [this is the case in which the computation
  of $\R f_*f^\times C$ will carry useful information about $f^\times C$],
what is relevant for us is that
$\R f_*f^\times C\cong\RHHom_{\Dqc(Y)^{}}^{}(\R f_*\co_X^{},C)$ is given
by the more computable expression if one of the conditions below
holds:
\be
\item
$f$ is finite and $C$ is bounded below.
\item
$f$ is finite and of finite Tor-dimension, and $C$ is arbitrary.
\ee
In passing we note that, historically, the most widely used special case of the above has
been where $f$ is a closed immersion---possibly of finite Tor-dimension.
\ermk

There isn't a whole lot more concrete computational knowledge about $f^\times$:
in the body of the article we
told the reader what is known when $f$ is smooth and proper, and
Remark~\ref{R95.3.1} tells us some more when $f$ is a finite map, possibly
of finite Tor-dimension. Nevertheless one can use these
(admittedly limited) pieces of information to deduce
useful facts. But for this it helps to know another functor,
a close cousin of $f^\times$. Before
we introduce it, a reminder might help.

In Reduction~\ref{R30.37} we met the following situation:
we were given composable
morphisms of schemes $U\stackrel u\la X\stackrel f\la Y$, with $u$ an open
immersion and $f$ proper. The Reduction was all about
computing the functor $\LL u^*f^\times$, and until now we haven't
mentioned that $\LL u^*f^\times$ depends only on the composite map
$U\stackrel{fu}\la Y$. That is:

\rmd{R95.-77}
If $g:U\la Y$ is any separated morphism
of finite type, Nagata's theorem~\cite{Nagata62,Conrad07} allows
us to choose a factorization of $g$ as
$U\stackrel u\la X\stackrel f\la Y$, with $u$ an open
immersion and $f$ proper. We then define $g^!=\LL u^*f^\times$, and it
is a theorem that $g^!$ is independent of the factorization up to
canonical isomorphism\footnote{If one follows the Grothendieck-Hartshorne
path to the subject, which we will discuss a little more in 
Remark~\ref{R95.-1}, then the way to see the isomorphism 
$g^!\cong\LL u^*f^\times$
is that the theory sets up a functor $g^!$, one shows
that for open immersions $u$ there is a natural isomorphism
$u^!\cong \LL u^*$, for proper maps $f$
there is a natural isomorphism $f^!\cong f^\times$, and
for composable maps there is a natural 
isomorphism $g^!=(fu)^!\cong u^!f^!$. But both
Deligne~\cite{Deligne66} and Verdier~\cite{Verdier68} sketched an argument
hinting how to prove directly the independence of factorization
of $\LL u^* f^\times$, using only~\ref{SS30.5.1} applied to suitable 
cartesian squares. There is 
some more detail in Lipman's book~\cite{Lipman09}, and for 
a fuller argument, which works for algebraic stacks 
and hence must handle the 2-category technicalities more carefully, the reader
is referred to~\cite{Neeman13}.}. Classically the bulk of Grothendieck
duality has been devoted to the study of the functors $g^!$.
See Remark~\ref{R95.-1} below.

Having introduced $g^!$, we are now in a position to
restate the recent progress in terms of it. Suppose
we are given, as in the last paragraph, a map of
schemes $g:U\la Y$ which factorizes as
$U\stackrel u\la X\stackrel f\la Y$, with $u$ an open
immersion and $f$ proper. The square
\[\xymatrix{
U \ar[r]^\id \ar[d]_\id & U \ar[d]^u\\
U\ar[r]^-u & X
}\]
is cartesian, hence it has a base-change map
$\Phi:u^\times\la\LL u^*$ as in Construction~\ref{C30.5}. The recent discoveries
can be summarized as saying
\be
\item
  The induced map
  $\Phi f^\times:u^\times f^\times\la \LL u^* f^\times$ is independent
of the factorization, and gives an unambiguous map
$g^\times\la g^!$. The construction taking $g$ to $g^!$
yields a 2-functor we will denote $(-)^!$---in other words there is
a compatibility with composition. And the map $g^\times\la g^!$
is compatible too, it is a morphism of
2-functors $\psi:(-)^\times\la(-)^!$ with many reasonable naturality
properties. The reader can find this
worked out in~\cite{Iyengar-Lipman-Neeman13},
and in greater generality in~\cite{Neeman13}.
\item
There are interesting situations in which some natural
functor $\Gamma$ takes $\psi(g):g^\times\la g^!$ to an isomorphism.
We have encountered two examples, 
namely \ref{ST30.39.1} and Sketch~\ref{SK35.7}.
For the general theory see~\cite{Iyengar-Lipman-Neeman13,Neeman13}.
\ee
\ermd

\apl{A95.-1024}
Suppose $f:X\la Y$ is proper.
Using the 2-functor $(-)^!$ one can prove the following 
\be
\item
If $G\in\dcohp(Y)$ then $f^\times G\in\dcohp(X)$. Here
$\dcohp\subset\Dqc$ is the full subcategory of all objects
with coherent cohomology sheaves, which vanish in sufficiently negative
degrees.
\setcounter{enumiv}{\value{enumi}}
\ee
Now suppose $f$ is not only proper, but also of finite Tor-dimension.
Then 
\be
\setcounter{enumi}{\value{enumiv}}
\item
If $G\in\dcoh(Y)$ then $f^\times G\in\dcoh(X)$. Here
$\dcoh\subset\Dqc$ is the full subcategory of complexes with
bounded, coherent cohomology sheaves.
\item
$G\in\dcohm(Y)$ then $f^\times G\in\dcohm(X)$. Here
$\dcohm\subset\Dqc$ is the full subcategory of complexes with
coherent cohomology sheaves, vanishing in sufficiently positive degrees.
\ee
Thus if $f$ is proper the adjoint pair
$\xymatrix{\R f_*:\Dqc(X)\ar@<0.5ex>[r] & \ar@<0.5ex>[l] \Dqc(Y):f^\times}$
restricts to an adjunction
$\xymatrix{\R f_*:\dcohp(X)\ar@<0.5ex>[r] & \ar@<0.5ex>[l] \dcohp(Y):f^\times}$.
When $f$ is not only proper but also of finite Tor-dimension, we also
have two more adjoint pairs, namely
$\xymatrix{\R f_*:\dcohm(X)\ar@<0.5ex>[r] & \ar@<0.5ex>[l] \dcohm(Y):f^\times}$
as well as
$\xymatrix{\R f_*:\dcoh(X)\ar@<0.5ex>[r] & \ar@<0.5ex>[l] \dcoh(Y):f^\times}$.

And the proofs of (i), (ii) and (iii) go as follows: the assertions are
local in $X$, hence it suffices to show that, for $u:U\la X$ an
open immersion from a (sufficiently
small) open affine subset, the map $\LL u^*f^\times=(fu)^!$
satisfies the properties. But $(-)^!$ is a 2-functor. We are allowed to factor
$fu$ in some other way, for example as a composite $fu=ghk$, and if
we choose our factorization wisely then $(fu)^!=k^!h^!g^!$ might be more
computable. For example $g,h,k$ might fall into the classes
where we understand $(-)^!$: open immersions, maps that are smooth and proper,
and finite morphisms (possibly of finite Tor-dimension). For details
the reader is
referred to \cite[Lemma~3.12 and its proof]{NeemanTIFR}.

Thus the abstract nonsense approach does recover ``coherent duality'' as
it was traditionally understood---meaning about complexes with coherent
cohomology.
The reader might also wonder about the relation of
(for example) the categories $\dcoh(X)$ and 
$\D^b(\mathbf{coh}/X)$. The relation is well-understood by now,
but falls outside the scope of Grothendieck duality---it isn't about the
functors $f^\times$ or $f^!$, it's a formal question about the
interplay among the myriad derived categories one can associate
to a single scheme $X$. We leave this out of
the survey---we've barely started Section~5, and we've already been
bombarded with a hail of new categories and functors. 
\eapl

\rmk{R95.-1}
Now that we have introduced $(-)^!$ we can recall a
historical point.
The Grothendieck approach is entirely in terms of the 2-functor $(-)^!$, the
2-functor $(-)^\times$ never comes up. That is: starting with a separated,
finite-type morphism
of noetherian schemes $f:X\la Y$, Grothendieck goes through an
intricate, arduous procedure
to arrive at a functor which turns out to agree with the $f^!$ of
Reminder~\ref{R95.-77}.
In the Grothendieck setting it's a major theorem that, for proper $f$,
the functor
$f^!$ is right adjoint to $\R f_*$, that is $f^!$ satisfies
the defining property
of $f^\times$. Of course if you take Remark~\ref{R95.-77} as the definition
of $f^!$ then this theorem becomes trivial: for a proper
$f$ we choose the Nagata factorization
$X\stackrel\id\la X\stackrel f\la Y$, and by definition
$f^!=\id^*f^\times=f^\times$.

We have already mentioned that, for general $f$, the 
existence of the right adjoint $f^\times$ for the functor
$\R f_*$  was 
first proved in Deligne~\cite{Deligne66}.
The article \cite{Deligne66} notes that $f^\times$ agrees with $f^!$ for
proper $f$, and for non-proper $f$ Deligne dismisses
$f^\times$ as too undeserving even to be graced with a
name, a functor that doesn't lend
itself to calculation.
For non-proper $f$
the functor $f^\times$ was deemed worthless and consigned to the trash
heap of obscurity, until the 1990s it remained nameless.

Until the 1980s, the
Deligne-Verdier approach to Grothendieck duality
existed only as an offhand aside
in Verdier~\cite{Verdier68}, a remark saying such a theory should be possible.
No one took the trouble to check the details.
Lipman plunged into this project sometime in the late 1980s---it took
him the best part of two decades, and the outcome
was the book~\cite{Lipman09}. In his development
of the subject Lipman chose not to ignore the right adjoint of $\R f_*$
when $f$ isn't proper. He brushed off Deligne's disparaging appraisal
of this functor---he
christened\footnote{The drawback of Lipman's notation
  is that, when handwritten, the symbol $f^\times$ is barely distinguishable
  from $f^*$. This renders it a calligraphic
  challenge to give blackboard talks in the subject.} it 
$f^\times$, and then went on 
to study its properties.

As it turns out Lipman was right in resurrecting the maligned
$f^\times$ from oblivion:
the key to the new progress is to
study $f^\times$ for non-proper maps. In
particular the computability of $f^\times$, when $f$ is a morphism of
affine schemes, plays a key role---see Remark~\ref{R30.40}.
\ermk

\rmk{R95.1}
Although none of the results presented in this survey is new, the
published expositions
of the recent work are all addressed to the expert---they are all couched in
terms of $g^!$, its relation with $g^\times$, and how to use the relation for
computations. Thus, even though the results give new and much
simplified proofs of concrete, old theorems,
you have to be an expert to discern this from the available
manuscripts. This was inevitable: after all
when we write research papers we'd like to have
them published. When the subject happens to be out of fashion this adds
an extra hurdle---the editor will naturally worry that the paper is
likely to have only a miniscule readership and a meager impact (as measured by
citations).
But, even with the most sympathetic
of editors, the paper will go to referees who will undoubtedly be experts.
Thus the authors will invariably try to write to impress the experts.
The focus will be on new results, not on simple proofs of old
theorems---and even when such simplifications
are present the fact may well be hidden, buried deep under a
mountain of
technicalities.

That said, the current survey is an attempt to open up the field to
non-specialists. If the subject is to have a revival then it's imperative
for the main points to be widely accessible. To keep this document
as readable as possible we have, up to Section~\ref{S95},
avoided all but the most indispensable machinery. We
are about to delve into areas where progress now seems within reach,
and this is forcing us
to first recall a little more background.
We have to catch up a little, before the open questions can be stated
clearly.

We have explained the 2-functor $(-)^!$. The next step is
\ermk

\rmd{R95.3}
We begin with a couple of old definitions
\be
\item
Let $X$ be a noetherian scheme. An object $D\in\dcoh(X)$ is a \emph{dualizing
  complex} if $\RHHom(-,D)$ takes $\dcoh(X)$ to
itself\footnote{In Remark~\ref{R95.3.1} we noted that there are
  two classical functors $\RHHom$, namely the one with values in
  $\Dqc(X)$ and the one with values in $\D(X)$. Since we are
only considering $\RHHom(C,D)$ with $C\in\dcoh(X)$ and $D\in\dcoh(X)\subset\Dqcpl(X)$, we are in a situation where the two are classically known to agree.}, and furthermore
the natural map $\co_X^{}\la\RHHom(D,D)$ is an isomorphism.
Equivalently: the functor $\RHHom(-,D)$ induces
an equivalence $\dcoh(X)\op\la\dcoh(X)$.
\item
Let $f:X\la Y$ be a flat, finite-type morphism of noetherian schemes,
and consider the diagram
\[\xymatrix@C+10pt@R+0pt{
  X\ar[r]^-\delta 
& X\times_Y^{} X \ar[r]^-{\pi_1^{}} \ar[d]_{\pi_2^{}} & X\ar[d]^f \\
& X\ar[r]^-f & Y  
}\]
A dualizing complex
$D\in\dcoh(X)$ is \emph{$Y$--rigid} if it comes together with an
isomorphism $D\la\delta^\times\big[\LL\pi_1^*D\oo^\LL_{}\LL\pi_2^*D\big]$.
\setcounter{enumiv}{\value{enumi}}
\ee 
An old theorem tells us
\be
\setcounter{enumi}{\value{enumiv}}
\item
Assume $X\stackrel f\la Y\stackrel g\la Z$ are finite-type,
flat, separated morphisms of noetherian
schemes. Then $f^!:\Dqc(Y)\la\Dqc(X)$ takes $Z$--rigid dualizing
complexes to $Z$--rigid dualizing complexes.
\ee
For the sake of historical accuracy:
the fact that $f^!$ takes dualizing complexes to dualizing
complexes goes all the way back to the dawn of the theory. The concept of
rigid dualizing complexes started with Van den Bergh~\cite{vandenBergh97}.
The formulation given here follows Lipman's reworking of Van den Bergh's
result. And Yekutieli and Zhang~\cite{Yekutieli-Zhang05,
  Yekutieli-Zhang06,Yekutieli-Zhang08,Yekutieli-Zhang09}
pursued this in depth, it was their contribution to simplifying
and extending the Grothendieck approach to the subject.

In other words: in an alternative universe, a survey of the field would
begin with rigid dualizing complexes and build up from there. It just so
happens that the Deligne-Verdier angle on the subject was the first
to achieve the status of satisfactorily cleaning up the foundations.
\ermd

\rmd{R95.5}
We should also remind the reader that there is an
interesting noncommutative version.
Since I'm not quite sure what a general noncommutative scheme should be,
I will confine the discussion to affine noncommutative schemes.

Let $R$ be a noetherian, commutative ring and let $S$ be a flat,
finitely generated, associative $R$--algebra which is
right and left noetherian (but not necessarily
commutative). Set $\se=S\oo_R^{}S\op$.
Following Yekutieli~\cite{Yekutieli92}, a dualizing complex
is an object $D\in\D(\se)$ such that the functor
$\RHom(-,D)$ yields an equivalence
of categories $\D^b(\mmod S)\op\la\D^b(\mmod {S\op})$.
The paper \cite{Yekutieli92} goes on to study the graded situation
and produce some examples of dualizing complexes.
If $S$ is commutative then dualizing complexes in $\D(\se)$,
in Yekutieli's sense,
can be shown to agree
with dualizing complexes in $\dcoh\big(\spec S\big)$,
as recalled in
Reminder~\ref{R95.3}(i).

Following  Van den
Bergh~\cite{vandenBergh97},
the dualizing complex $D\in\D(\se)$ is $R$--\emph{rigid} if it comes
together with
an isomorphism
\[
D\quad\cong\quad \RHom_\se^{}\big(S\,\,,\,\,D\oo_R^\LL D\big)\ .
\]
If $R$ is a field, and $S$ has a filtration whose associated
graded ring is commutative and finitely generated as an $R$--algebra, then
\cite{vandenBergh97} cleverly shows that a rigid dualizing complex
exists.
When $S$ is commutative then an $R$--rigid dualizing
complexes in $\D(\se)$, in Van den Bergh's sense,
can be shown to agree with the
$\spec R$--rigid dualizing complex in $\dcoh\big(\spec S\big)$,
as
recalled in
Reminder~\ref{R95.3}(ii).

There has been literature pursuing this further, for a
couple of early papers the reader is referred to 
Yekutieli and Zhang~\cite{Yekutieli-Zhang97,Yekutieli-Zhang99}. 
But it's now high time to move on to the open problems.
\ermd

\subsection{Foundational questions}
\label{SS95.2}

And now we are ready to state the first open question. This one is
based not on the very recent work, it derives from the formulas
of Avramov and Iyengar~\cite{Avramov-Iyengar08,Avramov-Iyengar-Lipman-Nayak10}
that finally opened our eyes. Fittingly the first open question is
about noncommutative algebraic geometry, which is unquestionably a hot
field nowadays.

\plm{P95.7}
Let $g:k\la R$ be any finite-type,
flat homomorphism of noetherian, commutative rings.
Let the relation between the
rings $R$ and $S$ be as in Reminder~\ref{R95.5}, meaning 
$S$ is an $R$--algebra satisfying all the hypotheses of Reminder~\ref{R95.5}. 
Let $N\in\D(\mmod R)\cong\dcoh\big(\spec R\big)$ be a $k$--rigid
dualizing complex.

Question: is $S\oo_\se^\LL\RHom_R^{}(S,S\oo_R^\LL N)$ a $k$--rigid dualizing
complex in $\D(\mmod \se)$?

We should remark that, if $S$ is commutative, this follows from
the isomophism $f^!N\cong S\oo_\se^\LL\RHom_R^{}(S,S\oo_R^\LL N)$ of
the Reduction Formula, which we met in~\S\ref{SS35.5}, coupled with 
Reminder~\ref{R95.3}(iii). But we don't yet understand the theory
well enough to have a simple, direct proof of
Reminder~\ref{R95.3}(iii) in the affine (commutative) case, hence have no idea
if the result can be extended to the noncommutative context.

We should note that the case $k=R$ is already interesting. In fact let
us confine ourselves to the case where the ring $k=R$ is Gorenstein; in this
case it is known that
$R\in\D(\mmod R)$ is an $R$--rigid dualizing complex, and the
question specializes to: is $S\oo_\se^\LL\RHom_R^{}(S,S)$ an $R$--rigid dualizing
complex in $\D(\mmod\se)$? The reader should note that,
in the noncommutative setting, $R$--rigid dualizing complexes
are known to exist only when $R$ is a field,
in particular the existence results to date
all assume equal characteristic.
\eplm

\plm{P95.999}
The work of Avramov, Iyengar and
Lipman~\cite[Section 3]{Avramov-Iyengar-Lipman11} suggests an alternative
definition for rigid dualizing complexes.
With $k\la R\la S$ as in Problem~\ref{P95.7}, we can declare that a dualizing
complex $D\in\D(\se)$ is \emph{AIL--$k$--rigid} if it comes with an isomorphism
\[
D\quad\cong\quad\RHom_S^{}\Big[\RHom_S^{}\big(D\,,\,S\oo_\se^\LL
  \RHom_k^{}(S,S)\big)\,,\,D\Big]
\]
From the Reduction Formula of Avramov, Iyengar,
Lipman and Nayak~\cite[(4.1.1)]{Avramov-Iyengar-Lipman-Nayak10}
we have that, as long as $k$ is regular and finite dimensional and
$S$ is commutative, the two notions of rigidity
agree.

Question 1: Do the two notions coincide when $S$ isn't commutative?

Question 2, assuming the notions are different: Suppose $N\in\D(R)$
is an AIL--$k$--rigid dualizing complex.
Is $S\oo_\se^\LL\RHom_R^{}(S,S\oo_R^\LL N)$ also an AIL--$k$--rigid dualizing
complex in $\D(\mmod \se)$?
\eplm

The third foundational problem is about a derived stack version of the
theory---as it happens derived stacks are also much in vogue nowadays.

\plm{P95.9}
In many of the theorems we had to assume flatness, or at the
very minimum finite
Tor-dimension. The modern way to get around this is to work in
the setting of derived algebraic geometry.

Question: is there an incarnation of the theory in derived algebraic
geometry?

The
Yekutieli school was the first to successfully employ DG methods in
Grothendieck duality: see Yekutieli and
Zhang~\cite{Yekutieli-Zhang04,
  Yekutieli-Zhang08,Yekutieli-Zhang09},
Yekutieli~\cite[survey]{Yekutieli10}, and more
recently Shaul~\cite{Shaul17}. The Lipman school, inspired by the
successes of the Yekutieli school, followed suit:
the affine case of the Reduction Formulas,
of Avramov and Iyengar, does extend from the flat case presented
in~\S\ref{SS35.5} to the case where the map $R\la S$ is of finite
Tor-dimension. The proof given in
Avramov, Iyengar, Lipman and
Nayak~\cite[Section~4]{Avramov-Iyengar-Lipman-Nayak10} goes by way of
differential graded algebras. See
also~\cite{Avramov-Iyengar-Lipman11} for further instances, of 
the Lipman school exploiting the DG methods
introduced by Yekutieli and Zhang. 

For some planned future projects see 
Yekutieli~\cite{Yekutieli19A}. Lipman is also interested in
pursuing further
the
methods of derived
algebraic geometry---he has been working his way through
Lurie's book~\cite{Lurie09}---but I'm not
aware of any manuscripts yet. In any case: at this point
the subject is in its infancy, all are welcome to join in.
\eplm

\rmk{R95.11}
The work in Hafiz Khusyairi 2017 PhD thesis might be relevant to
Problem~\ref{P95.9}---the thesis is
entirely in the setting of old-fashioned, ordinary, commutative
schemes, but uses the formulas of Avramov and
Iyengar as the starting point for setting up the theory. It then proceeds
to develop the usual functoriality properties from there. Since the formulas
have a DG analog Khusyairi's work might generalize.
\ermk

\subsection{Computational problems}
\label{SS95.3}

In the previous section we sketched three foundational problems, about
extending the theory to noncommutative and to derived algebraic geometry---both
of which are ``in'' fields nowadays.

Let us now return to more classical problems. In the old, traditional
world of ordinary, commutative algebraic geometry the foundations of
Grothendieck duality have reached a reasonably satisfactory state.
It is feasible to introduce the players and describe the relations
among them in what could plausibly be
called a short space, and it is possible
to do so in such a way that the traditional computations become
transparent and brief.

But the problem is that the traditional computations are limited. Let us
assume $f:X\la Y$ proper and
of finite Tor-dimension, in which case Remark~\ref{R-1.7} allows us
to reduce the problem to computing $f^\times\co_Y^{}$ and the
map $\e:\R f_*f^\times\co_Y^{}\la\co_Y^{}$. The classical literature
gives us a good understanding in the case where $f$ is smooth,
and an understanding of some sort in situations that are easily reduced
to the smooth case, for example when $f$ is Cohen-Macauley. Beyond
that, what's known is not all that useful.

As it happens Nayak and Sastry are in the process of
writing up a comprehensive account of what is known about the
computations. Before long there will be
a manuscript containing everything that has
been figured out so far---it will be valuable
to have it all assembled in one place and the connections worked out.
Once the document is ready the interested
reader will be able to see, in print and in detail,
just how paltry our understanding
really is.

\rmk{R95.3.99}
The last paragraphs should not be interpreted as belittling the
traditional 
case of
Grothendieck duality, the special case where
$f:X\la Y$ is smooth and proper---this classical
situation is already fascinating and
has spawned a rich literature spanning
many decades. Specializing further, assume that
$Y=\spec k$ is a point
and $f:X\la Y$ is smooth, proper
and of relative dimension 1, and we find ourselves in ancient
territory---we're in the
well-understood case of duality for curves. The reader
can find an excellent
exposition of the old approaches in
Serre~\cite[pp.~25-34 and pp.~76-81]{Serre59}, and
[as far as the author knows] the cleverest, most
recent 
idea is already five decades old, see
Tate~\cite{Tate68}.
The case where $Y=\spec k$ is still a point, $f:X\la Y$ is still smooth and
proper, but the relative dimension is arbitrary is classical Serre duality,
the reader is referred to Serre~\cite{Serre55A}, and also to the sketch
presented in Section~\ref{SS29.987}. In Be{\u\i}linson~\cite{Beilinson80}
we learn how
to generalize Tate's clever trick to higher dimension.

If $k=\cpl$ then $X$ is a smooth, compact K\"ahler manifold, and the
Hodge decomposition theorem identifies $H^n(\Omega^n_f)=\R f_*\Omega_f^n[n]$ with
$H^{n,n}(X)\cong H^{2n}(X,\cpl)\cong\cpl$, where the last isomorphism
is by Poincar\'e duality. We would therefore expect to be able to
understand the residue map from this perspective too. The reader can
find this explored in Harvey~\cite{Harvey70}, Tong~\cite{Tong73},
and more recently in Sastry and Tong~\cite{Sastry-Tong03}.

Now let us return to the generality of the relative case: that is
$f:X\la Y$ is assumed smooth and proper but $Y$
is an arbitrary noetherian scheme. We know, from the results surveyed in
this article, that $f^\times\co_Y^{}$ is canonically isomorphic
to $\Omega^n_{f}[n]$, and that the counit of adjunction
$\e:\R f_*f^\times\co_Y^{}\la\co_Y^{}$ is determined by the map
taking a relative meromorphic $n$-form to its residue. In this
article we presented a very recent
approach to these theorems, we should
say something about the older methods---after all understanding the
relationship of the old tack with the new might well prove fruitful
and illuminating.

The first issue is that all the data must be compatible with
composition. That is: if $X\stackrel f\la Y\stackrel g\la Z$ are
composable morphisms of schemes, both of which are smooth and proper,
then the composite is smooth and proper and we have a string of
canonical
isomorphisms
\[\begin{array}{ccccl}
\Omega^{m+n}_{gf}[m+n]&\cong&(gf)^\times\co_Z^{}&\cong &f^\times g^\times\co_Z^{}\\
&&&\cong&
\LL f^*g^\times\co_Z^{}\oo^\LL_{}f^\times\co_Y^{}\\
&&&\cong&\LL f^*\Omega^m_{g}[m]\oo^\LL\Omega^n_{f}[n]
\end{array}\]
The reader might wonder whether the composite is the
obvious isomorphism---not sursprisingly the answer turns out to be Yes,
see Lipman and Sastry~\cite{Lipman-Sastry92}.
Furthermore the counits of adjunction must be compatible. We have a
counit of adjunction $\e(gf):\R(gf)_*(gf)^\times\co_Z^{}\la\co_Z^{}$, which must
agree with the composite
\[
\xymatrix@C+30pt{
\R g_*\R f_*f^\times g^\times\co_Z^{}\ar[r]^-{\e(f)} &
\R g_*g^\times\co_Z^{}\ar[r]^-{\e(g)} &\co_Z^{}\ .
}\]
Rewriting this in terms the string of canonical isomorphisms above yields
a diagram which must commute
\[
\xymatrix@C+20pt{
\R g_*\big(\Omega^m_{g}[m]\oo^\LL\R f_*\Omega^n_{f}[n]\big)\ar[rr]^-{\cong}  
\ar[d]_-{\e(f)} &&
\R(gf)_*\Omega^{m+n}_{gf}[m+n]
\ar[d]^-{\e(gf)} \\
\R g_*\big(\Omega^m_{g}[m]\oo^\LL\co_Y^{}\big) \ar@{=}[r] &
\R g_*\Omega^m_{g}[m]
\ar[r]^-{\e(g)} &
\co_Z^{}
}\]
and the commutativity can be interpreted as a compatibility condition
on the residue maps.

In the previous paragraph we learned that the counit of adjunction
$\e:\R f_*\Omega^n_{f}[n]\la\co_Y^{}$, and hence the closely related
residue map $\rho:\R f_*\Gamma_W^{}\Omega^n_{f}[n]\la\co_Y^{}$,
must be compatible with composition. 
It's even easier to see that $\rho$ must be compatible with
flat base change. For the purpose of computations,
the compatibility with flat base change allows
us to assume that $Y=\spec R$ is affine---and if it helps we may even assume
that $R$ is a
(strictly) henselian or even a complete local ring. Grothendieck's
GFGA result~\cite[Th\'eor\`eme~5.1.4]{Grothendieck61} allows us to 
replace $X$ by its formal completion, and for some time now
the experts
have been pursuing the idea that doing so might lead to  
a better understanding of $\rho$.
For a more extensive treatment the reader is referred to
the forthcoming article by Nayak and Sastry;
but see also
Alonso, Jerem{\'{\i}}as and
Lipman~\cite{Alonso-Jeremias-Lipman97,Alonso-Jeremias-Lipman99},
Lipman Nayak and Sastry~\cite{Lipman-Nayak-Sastry05},
Nayak~\cite{Nayak05}, Nayak and Sastry~\cite{Nayak-Sastry17} and
Sastry~\cite{Sastry05}.
\ermk

\rmk{R95.3.109}
In  the
special case of smooth and proper morphisms $f:X\la Y$,
Remark~\ref{R95.3.99} surveyed some of the work 
done in the quest for a better understanding
of Grothendieck duality.  The
opening paragraphs of Section~\ref{SS95.3} dismissed what's
known about more general $f$ as being of limited computational value.

There is some literature: the reader might wish to look at
Huang~\cite{Huang00,Huang01},
Kersken~\cite{Kersken83}, Par{\v s}in~\cite{Parsin76} and 
Yekutieli~\cite{Yekutieli92A} (see also the appendix by Sastry).
\ermk

\plm{P95.3.119}
Now put the recent results at center stage---they should allow
us to go further with the computations.
At least when $f$ is flat,
we have a simple and explicit formula for $f^!$.
The generalization of Reduction~\ref{R30.37} gives an isomorphism
$f^!=\LL\delta^*\pi^\times\LL f^*$,
where $\delta:X\la X\times_Y^{}X$ is the diagonal map
and $\pi:X\times_Y^{}X\la X$ is the (second) projection.
Any colocalization $c:\Gamma\la\id$,
where $\Gamma$
takes the map $\psi:f^\times\la f^!$ to an isomorphism, will permit us
to form the composite
\[\xymatrix@C+30pt{
\R f_*\Gamma f^! \ar[r]^-{\R f_*(\Gamma\psi)^{-1}} & \R f_*\Gamma f^\times
\ar[r]^-{\R f_* cf^\times} & \R f_*f^\times\ar[r]^-\e &\id
}\]
which should be computable, at least in the special case where
$X$ and $Y$ are affine. For suitable
choices, of the colocalization $c:\Gamma\la\id$, the
composite should deliver useful information about $\e$---and 
those of us competent to carry out the computations
should be able to learn much more about the map
$\e:\R f_*f^\times\la\id$.

The computations will involve Hochschild
homology and cohomology---terms like $S\oo_\se^{}\RHom_R^{}\big(S,S\oo_R^{}N\big)$
are bound to appear. Fortunately the world is full of experts in
Hochschild homology and cohomology, and once they take an
interest they will undoubtedly be able to move these computations much
further than the handful of us, the few people who have been working on
Grothendieck duality. Let's face it: in our tiny group none is adept at
handling the Hochschild machinery. The Hochschild
experts should feel invited to move
right in.
\eplm

\appendix

\section{A computation of the base-change map $u^\times\la\LL u^*$ when
$u:U\la C$ is an open immersion of curves}
\label{A97}

Let $k$ be a field, let $C$ be a
complete algebraic curve smooth over $k$, and let
$p\in C$ be a $k$--rational point. Put $U=C-\{p\}$ and let
$u:U\la C$ be the open immersion. The square
\[\xymatrix{
U\ar@{=}[r]\ar@{=}[d] & U\ar[d]^u \\
U\ar[r]^-u & C
}\]
is cartesian and the horizontal maps are flat, and Construction~\ref{C30.5}
yields a base-change map $\Phi:u^\times\la \LL u^*$.
Let $\cl$ be a line bundle on
$C$; we propose
to compute the map $\Phi(\cl):u^\times\cl\la \LL u^*\cl$.
In Remark~\ref{R95.3.1} we met the isomorphism
\[
\R u_*u^\times \cl\quad\cong\quad
\RHHom_{\Dqc(C)^{}}^{}(\R u_*\co_U^{},\cl)
\]
One may check that the counit of adjunction $\e:\R u_*u^\times \cl\la \cl$
is the map obtained by applying the functor $\RHHom_{\Dqc(C)^{}}^{}(-,\cl)$
to the morphism $\co_C^{}\la\R u_*\co_U^{}$.
And the map $u^\times\la\LL u^*$ is just the composite
$u^\times\stackrel{(\e')^{-1}u^\times}\la
\LL u^*\R u_*u^\times \stackrel{\LL u^*\e}\la \LL u^*$,
where $\e':\LL u^*\R u_*\la\id$ is the (invertible) counit of adjunction.

So much for abstract nonsense. Concretely we are reduced 
to computing what the functor $\RHHom_{\Dqc(C)^{}}^{}(-,\cl)$
does to the morphism $\co_C^{}\la\R u_*\co_U^{}$, after which we
will apply
$\LL u^*$. Now the map $\co_C^{}\la\R u_*\co_U^{}=u_*\co_U^{}$
is the direct limit, as $n\sr\infty$, of the maps $\co_C^{}\la\co_C^{}(np)$.
This means that, in the derived category $\Dqc(C)$, we need to
compute the homotopy inverse limit of the sequence $\cl(-np)\la\cl$.
This is what we will now do.

Let $R=\co_{C,p}^{}$, that is the stalk at $p$ of the structure sheaf
$\co_C^{}$. Let $\m\subset R$ be the maximal ideal. For
each $n$ we have a triangle
\[\xymatrix{
\cl(-np)\ar[r] &\cl\ar[r] & \cl\oo R/\m^n
}\]
and taking homotopy inverse limits over $n$ will yield a triangle.
In general I find homotopy inverse limits difficult, but
in the case of the inverse system $\cl\oo R/\m^n$ it isn't so
bad.

Let us first take the homotopy inverse limit in the category $\D(C)$,
where we allow all complexes of sheaves of $\co_C^{}$--modules---not
only ones with quasicoherent cohomology, see Remark~\ref{R95.3.1}
for a discussion. Let $i:p\la C$
be the inclusion of $p$; we turn it into a map of ringed spaces
by giving $p$ the structure sheaf $R$. The functor $i_*$ (extension by zero)
is exact,
and has an exact left adjoint---the functor taking a sheaf to its
stalk at $p$. Hence the induced functor $i_*:\D(p)\la\D(C)$ respects
products and therefore homotopy inverse limits. Thus
the homotopy inverse limit of the system
$\cl\oo R/\m^n$ can be computed in $\D(p)$, and it comes
down to the sheaf $\cl\oo i_*\wh R$, the extension by zero of the completion
of the stalk at $p$ of $\cl$.

This sheaf is manifestly not quasicoherent---to compute the homotopy
inverse limit in the category $\Dqc(C)$ we need to derived quasicoherate.
That is: we replace by an injective resolution and then quasicoherate.
The injective resolution is easy enough: if $K$ is the quotient field of
$R$ and $\wh K$ its $\m$-adic completion (i.e.~the quotient field of $\wh R$),
then an injective resolution of $\wh R$ as an $R$--module is given by
$\wh K\la \wh K/\wh R$, and $i_*\wh K\la i_*\big[\wh K/\wh R]$ is an
injective resolution of $i_*\wh R$ in the category of sheaves of
$\co_C^{}$--modules. The module $ i_*\big[\wh K/\wh R]$ is quasicoherent
as it stands, and the quasicoherator takes $i_*\wh K$ to the constant
sheaf $\wh K$. The morphism $\cl\la \holim (\cl\oo R/\m^n)$ becomes identified
with the cochain map
\[\xymatrix@C+20pt{
0\ar[r] & \cl\ar[r]\ar[d] & 0 \ar[r]\ar[d] & 0 \\
0 \ar[r]& \cl\oo\wh K \ar[r] & \cl\oo i_*\big[\wh K/\wh R] \ar[r] & 0
}\]
Applying the functor $\LL u^*=u^*$, that it restricting to $U\subset C$,
kills the sheaf $i_*\big[\wh K/\wh R]$. We deduce the map
of cochain complexes
\[\xymatrix@C+20pt{
0\ar[r] & u^*\cl\ar[r]\ar[d] & 0 \ar[r]\ar[d] & 0 \\
0 \ar[r]& u^*\cl\oo\wh K \ar[r] & 0 \ar[r] & 0
}\]
And the map $u^\times\cl\la \LL u^*\cl$ is obtained by
completing the triangle, it is the cochain map
\[\xymatrix@C+20pt{
0\ar[r] & u^*\cl\ar[r]\ar[d] &   u^*\cl\oo\wh K \ar[r]\ar[d] & 0 \\
0 \ar[r]& u^*\cl \ar[r] & 0 \ar[r] & 0
}\]
To sum it up: the curve $U$ is affine, let's say $U=\spec S$. The line
bundle $u^*\cl\in\Dqc(U)$ corresponds, under the equivalence
$\Dqc(U)\cong\D(S)$, to the rank--1 projective $S$--module
$L=\Gamma(U,\cl)$. Consider the short exact sequence
of $S$--modules
\[\xymatrix@C+20pt{
0\ar[r] & L\ar[r] &   L\oo\wh K \ar[r]&\frac{L\oo\wh K}{L} \ar[r] & 0
}\]
The  morphism $u^\times\cl\la\LL u^*\cl$ in the derived category
$\Dqc(U)$ corresponds, under the equivalence $\Dqc(U)\cong\D(S)$, to 
the map $\frac{L\oo\wh K}{L}[-1]\la L$ that is represented by
the short exact sequence.

Now there is an isomorphism of $S$--modules
$\Hom_k^{}(S,k)\cong\frac{L\oo\wh K}{L}$.
The bad way to see this is as follows: both are injective $S$--modules, and
the indecomposable injectives have the same multiplicity on both sides,
see~\cite{Neeman13B}. But in the case where $\cl$
is the canonical bundle $\Omega_C^1$ we know there is a canonical
isomorphism. If $f:C\la\spec k$ is the projection
to a point, then $f^\times k\cong \Omega_C^1[1]$ canonically, hence
\[
\Hom_k^{}(S,k)\quad\cong\quad 
(fu)^\times k\quad\cong\quad
u^\times f^\times k \quad\cong\quad
u^\times\Omega_C^1[1]
\quad\cong\quad \frac{\Omega_U^1\oo\wh K}{\Omega_U^1}
\]
where the isomorphisms are all canonical.

\def\cprime{$'$}
\providecommand{\bysame}{\leavevmode\hbox to3em{\hrulefill}\thinspace}
\providecommand{\MR}{\relax\ifhmode\unskip\space\fi MR }
\providecommand{\MRhref}[2]{%
  \href{http://www.ams.org/mathscinet-getitem?mr=#1}{#2}
}
\providecommand{\href}[2]{#2}


\begin{thebibliography}{10}

\bibitem{Alonso-Jeremias-Lipman97}
Leovigildo Alonso~Tarr{\'\i}o, Ana Jerem{\'\i}as~L{\'o}pez, and Joseph Lipman,
  \emph{Local homology and cohomology on schemes}, Ann. Sci. \'Ecole Norm. Sup.
  (4) \textbf{30} (1997), no.~1, 1--39.

\bibitem{Alonso-Jeremias-Lipman99}
\bysame, \emph{Studies in duality on {N}oetherian formal schemes and
  non-{N}oetherian ordinary schemes}, Contemporary Mathematics, vol. 244,
  American Mathematical Society, Providence, RI, 1999.

\bibitem{Alonso-Jeremias-Lipman11}
\bysame, \emph{Bivariance, {G}rothendieck duality and {H}ochschild homology
  {I}: {C}onstruction of a bivariant theory}, Asian J. Math. \textbf{15}
  (2011), no.~3, 451--497.

\bibitem{Alonso-Jeremias-Lipman12}
\bysame, \emph{Bivariance, {G}rothendieck duality and {H}ochschild homology,
  {II}: {T}he fundamental class of a flat scheme-map}, Adv. Math. \textbf{257}
  (2014), 365--461.

\bibitem{Alonso-Jeremias-Souto04}
Leovigildo Alonso~Tarr{\'\i}o, Ana Jerem{\'\i}as~L\'opez, and Mar{\'\i}a~Jos\'e
  Souto~Salorio, \emph{Bousfield localization on formal schemes}, J. Algebra
  \textbf{278} (2004), no.~2, 585--610.

\bibitem{Avramov-Iyengar-Lipman11}
Luchezar Avramov, Srikanth~B. Iyengar, and Joseph Lipman, \emph{Reflexivity and
  rigidity for complexes, {II}: {S}chemes}, Algebra Number Theory \textbf{5}
  (2011), no.~3, 379--429.

\bibitem{Avramov-Iyengar08}
Luchezar~L. Avramov and Srikanth~B. Iyengar, \emph{Gorenstein algebras and
  {H}ochschild cohomology}, Michigan Math. J. \textbf{57} (2008), 17--35,
  Special volume in honor of Melvin Hochster.

\bibitem{Avramov-Iyengar-Lipman-Nayak10}
Luchezar~L. Avramov, Srikanth~B. Iyengar, Joseph Lipman, and Suresh Nayak,
  \emph{Reduction of derived {H}ochschild functors over commutative algebras
  and schemes}, Adv. Math. \textbf{223} (2010), no.~2, 735--772.

\bibitem{Balmer-DellAmbrogio-Sanders15}
Paul Balmer, Ivo Dell'Ambrogio, and Beren Sanders, \emph{Grothendieck-{N}eeman
  duality and the {W}irthm\"uller isomorphism}, Compos. Math. \textbf{152}
  (2016), no.~8, 1740--1776.

\bibitem{Balmer-Favi11}
Paul Balmer and Giordano Favi, \emph{Generalized tensor idempotents and the
  telescope conjecture}, Proc. Lond. Math. Soc. (3) \textbf{102} (2011), no.~6,
  1161--1185.

\bibitem{Beilinson80}
Alexandre Be{\u\i}linson, \emph{Residues and ad\`eles}, Funktsional. Anal. i
  Prilozhen. \textbf{14} (1980), no.~1, 44--45.

\bibitem{Bokstedt-Neeman93}
Marcel B\"okstedt and Amnon Neeman, \emph{Homotopy limits in triangulated
  categories}, Compositio Math. \textbf{86} (1993), 209--234.

\bibitem{Bousfield79A}
A.K. Bousfield, \emph{The localization of spectra with respect to homology},
  Topology \textbf{18} (1979), 257--281.

\bibitem{Conrad00}
Brian Conrad, \emph{Grothendieck duality and base change}, Lecture Notes in
  Mathematics, vol. 1750, Springer-Verlag, Berlin, 2000.

\bibitem{Conrad07}
\bysame, \emph{Deligne's notes on {N}agata compactifications}, J. Ramanujan
  Math. Soc. \textbf{22} (2007), no.~3, 205--257.

\bibitem{Deligne66}
Pierre Deligne, \emph{Cohomology \`a support propre en construction du foncteur
  $f^!$}, Residues and Duality, Lecture Notes in Mathematics, vol.~20,
  Springer--Verlag, 1966, pp.~404--421.

\bibitem{Gabriel-Zisman67}
Peter Gabriel and Michel Zisman, \emph{Calculus of fractions and homotopy
  theory}, Ergebnisse der Mathematik und ihrer Grenzgebiete, Band 35,
  Springer-Verlag New York, Inc., New York, 1967.

\bibitem{Grothendieck61}
Alexandre Grothendieck, \emph{{\'El\'ements} de g\'eom\'etrie alg\'ebrique
  {III}. {{\'Etude}} cohomologique des faisceaux coh\'erents {I}}, Inst. Hautes
  \'Etudes Sci. Publ. Math. (1961), no.~11, 425--511.

\bibitem{Hartshorne66}
Robin Hartshorne, \emph{Residues and duality}, Lecture notes of a seminar on
  the work of A. Grothendieck, given at Harvard 1963/64. With an appendix by P.
  Deligne. Lecture Notes in Mathematics, No. 20, Springer-Verlag, Berlin, 1966.

\bibitem{Harvey70}
F.~Reese Harvey, \emph{Integral formulae connected by {D}olbeault's
  isomorphism}, Rice Univ. Studies \textbf{56} (1970), no.~2, 77--97 (1971).

\bibitem{Hochschild-Kostant-Rosenberg62}
Gerhard Hochschild, Bertram Kostant, and Alex Rosenberg, \emph{Differential
  forms on regular affine algebras}, Trans. Amer. Math. Soc. \textbf{102}
  (1962), 383--408.

\bibitem{Hopkins85}
Michael~J. Hopkins, \emph{Global methods in homotopy theory}, Homotopy
  Theory---Proceedings of the Durham Symposium 1985, London Math. Soc. Lecture
  Notes Series, vol. 117, Cambridge University Press, 1987, pp.~73--96.

\bibitem{Huang00}
I-Chiau Huang, \emph{An explicit construction of residual complexes}, J.
  Algebra \textbf{225} (2000), no.~2, 698--739.

\bibitem{Huang01}
\bysame, \emph{The residue theorem via an explicit construction of traces}, J.
  Algebra \textbf{245} (2001), no.~1, 310--354.

\bibitem{Hubl-Sastry93}
Reinhold H{\"u}bl and Pramathanath Sastry, \emph{Regular differential forms and
  relative duality}, Amer. J. Math. \textbf{115} (1993), no.~4, 749--787.

\bibitem{Illusie71C}
Luc Illusie, \emph{Conditions de finitude}, Th\'eorie des intersections et
  th\'eor\`eme de {R}iemann-{R}och, Springer-Verlag, Berlin, 1971, S\'eminaire
  de G\'eom\'etrie Alg\'ebrique du Bois-Marie 1966--1967 (SGA 6, Expos{\'e}
  III), pp.~222--273. Lecture Notes in Mathematics, Vol. 225.

\bibitem{Iyengar-Lipman-Neeman13}
Srikanth~B. Iyengar, Joseph Lipman, and Amnon Neeman, \emph{Relation between
  two twisted inverse image pseudofunctors in duality theory}, Compos. Math.
  \textbf{151} (2015), no.~4, 735--764.

\bibitem{Kersken83}
Masumi Kersken, \emph{Cousinkomplex und {N}ennersysteme}, Math. Z. \textbf{182}
  (1983), no.~3, 389--402.

\bibitem{Lipman84}
Joseph Lipman, \emph{Dualizing sheaves, differentials and residues on algebraic
  varieties}, Ast\'erisque (1984), no.~117, ii+138.

\bibitem{Lipman87}
\bysame, \emph{Residues and traces of differential forms via {H}ochschild
  homology}, Contemporary Mathematics, vol.~61, American Mathematical Society,
  Providence, RI, 1987.

\bibitem{Lipman02}
\bysame, \emph{Lectures on local cohomology and duality}, Local cohomology and
  its applications (Guanajuato, 1999), Lecture Notes in Pure and Appl. Math.,
  vol. 226, Dekker, New York, 2002, pp.~39--89.

\bibitem{Lipman09}
\bysame, \emph{Notes on derived functors and {G}rothendieck duality},
  Foundations of {G}rothendieck duality for diagrams of schemes, Lecture Notes
  in Mathematics, vol. 1960, Springer, Berlin, 2009, pp.~1--259.

\bibitem{Lipman-Nayak-Sastry05}
Joseph Lipman, Suresh Nayak, and Pramathanath Sastry, \emph{Pseudofunctorial
  behavior of {C}ousin complexes on formal schemes}, Variance and duality for
  Cousin complexes on formal schemes, Contemp. Math., vol. 375, Amer. Math.
  Soc., Providence, RI, 2005, pp.~3--133.

\bibitem{Lipman-Neeman07}
Joseph Lipman and Amnon Neeman, \emph{Quasi-perfect scheme maps and boundedness
  of the twisted inverse image functor}, Illinois J. Math. \textbf{51} (2007),
  209--236.

\bibitem{Lipman-Neeman18}
\bysame, \emph{On the fundamental class of an essentially smooth scheme-map},
  Algebr. Geom. \textbf{5} (2018), no.~2, 131--159.

\bibitem{Lipman-Sastry92}
Joseph Lipman and Pramathanath Sastry, \emph{Regular differential and
  equidimensional scheme-maps}, J. Algebraic Geom. \textbf{1} (1992), no.~1,
  101--130.

\bibitem{Lurie09}
Jacob Lurie, \emph{Higher topos theory}, Annals of Mathematics Studies, vol.
  170, Princeton University Press, Princeton, NJ, 2009.

\bibitem{Nagata62}
Masayoshi Nagata, \emph{Imbedding of an abstract variety in a complete
  variety}, J. Math. Kyoto Univ. \textbf{2} (1962), 1--10.

\bibitem{Nayak05}
Suresh Nayak, \emph{Pasting pseudofunctors}, Variance and duality for Cousin
  complexes on formal schemes, Contemp. Math., vol. 375, Amer. Math. Soc.,
  Providence, RI, 2005, pp.~195--271.

\bibitem{Nayak09}
\bysame, \emph{Compactification for essentially finite-type maps}, Adv. Math.
  \textbf{222} (2009), no.~2, 527--546.

\bibitem{Nayak-Sastry17}
Suresh Nayak and Pramathanath Sastry, \emph{Transitivity in duality for formal
  schemes i},  (2017), preprint.

\bibitem{Neeman13}
Amnon Neeman, \emph{An improvement on the base-change theorem and the functor
  $f^!$}, https://arxiv.org/abs/1406.7599.

\bibitem{NeemanTIFR}
\bysame, \emph{The relation between grothendieck duality and hochschild
  homology}, To appear in the proceedings of the TIFR international colloquium,
  2016.

\bibitem{Neeman92B}
\bysame, \emph{The chromatic tower for {$D(R)$}}, Topology \textbf{31} (1992),
  519--532.

\bibitem{Neeman96}
\bysame, \emph{The {Grothendieck} duality theorem via {Bousfield's} techniques
  and {Brown} representability}, Jour. Amer. Math. Soc. \textbf{9} (1996),
  205--236.

\bibitem{Neeman13B}
\bysame, \emph{The decomposition of {$\text{\rm Hom}_k(S,k)$} into
  indecomposable injectives}, Acta Math. Vietnam. \textbf{40} (2015), no.~2,
  331--338.

\bibitem{Neeman13A}
\bysame, \emph{Traces and residues}, Indiana Univ. Math. J. \textbf{64} (2015),
  no.~1, 217--229.

\bibitem{Parsin76}
A.~N. Par\v~sin, \emph{On the arithmetic of two-dimensional schemes. {I}.
  {D}istributions and residues}, Izv. Akad. Nauk SSSR Ser. Mat. \textbf{40}
  (1976), no.~4, 736--773, 949.

\bibitem{Porta-Shaul-Yekutieli14}
Marco Porta, Liran Shaul, and Amnon Yekutieli, \emph{On the homology of
  completion and torsion}, Algebr. Represent. Theory \textbf{17} (2014), no.~1,
  31--67.

\bibitem{Sastry04}
Pramathanath Sastry, \emph{Base change and {G}rothendieck duality for
  {C}ohen-{M}acaulay maps}, Compos. Math. \textbf{140} (2004), no.~3, 729--777.

\bibitem{Sastry05}
\bysame, \emph{Duality for {C}ousin complexes}, Variance and duality for Cousin
  complexes on formal schemes, Contemp. Math., vol. 375, Amer. Math. Soc.,
  Providence, RI, 2005, pp.~137--192.

\bibitem{Sastry-Tong03}
Pramathanath Sastry and Yue Lin~L. Tong, \emph{The {G}rothendieck trace and the
  de {R}ham integral}, Canad. Math. Bull. \textbf{46} (2003), no.~3, 429--440.

\bibitem{Serre55A}
Jean-Pierre Serre, \emph{Un th\'eor\`eme de dualit\'e}, Comment. Math. Helv.
  \textbf{29} (1955), 9--26.

\bibitem{Serre59}
\bysame, \emph{Groupes alg\'ebriques et corps de classes}, Publications de
  l'institut de math\'ematique de l'universit\'e de Nancago, VII. Hermann,
  Paris, 1959.

\bibitem{Shaul15}
Liran Shaul, \emph{Reduction of {H}ochschild cohomology over algebras finite
  over their center}, J. Pure Appl. Algebra \textbf{219} (2015), no.~10,
  4368--4377.

\bibitem{Shaul16}
\bysame, \emph{Relations between derived {H}ochschild functors via twisting},
  Comm. Algebra \textbf{44} (2016), no.~7, 2898--2907.

\bibitem{Shaul17}
\bysame, \emph{The twisted inverse image pseudofunctor over commutative {DG}
  rings and perfect base change}, Adv. Math. \textbf{320} (2017), 279--328.

\bibitem{Tate68}
John Tate, \emph{Residues of differentials on curves}, Ann. Sci. \'Ecole Norm.
  Sup. (4) \textbf{1} (1968), 149--159.

\bibitem{Thomason97}
Robert~W. Thomason, \emph{The classification of triangulated subcategories},
  Compositio Math. \textbf{105} (1997), 1--27.

\bibitem{ThomTro}
Robert~W. Thomason and Thomas~F. Trobaugh, \emph{Higher algebraic {K--theory}
  of schemes and of derived categories}, The Grothendieck Festschrift ( a
  collection of papers to honor Grothendieck's 60'th birthday), vol.~3,
  Birkh{\"a}user, 1990, pp.~247--435.

\bibitem{Tong73}
Yue Lin~L. Tong, \emph{Integral representation formulae and {G}rothendieck
  residue symbol}, Amer. J. Math. \textbf{95} (1973), 904--917.

\bibitem{vandenBergh97}
Michel Van~den Bergh, \emph{Existence theorems for dualizing complexes over
  non-commutative graded and filtered rings}, J. Algebra \textbf{195} (1997),
  no.~2, 662--679.

\bibitem{Verdier68}
Jean-Louis Verdier, \emph{Base change for twisted inverse images of coherent
  sheaves}, Algebraic {G}eometry ({I}nternat. {C}olloq., {T}ata {I}nst. {F}und.
  {R}es., {B}ombay, 1968), Oxford Univ. Press, London, 1969, pp.~393--408.

\bibitem{Yekutieli19A}
Amnon Yekutieli, \emph{The derived category of sheaves of commutative dg rings
  (preview)}, arXiv:1312.6411.

\bibitem{Yekutieli19}
\bysame, \emph{Duality and tilting for commutative dg rings}, arXiv:1312.6411.

\bibitem{Yekutieli92}
\bysame, \emph{Dualizing complexes over noncommutative graded algebras}, J.
  Algebra \textbf{153} (1992), no.~1, 41--84.

\bibitem{Yekutieli92A}
\bysame, \emph{An explicit construction of the {G}rothendieck residue complex},
  Ast\'erisque (1992), no.~208, 127, With an appendix by Pramathanath Sastry.

\bibitem{Yekutieli10}
\bysame, \emph{Rigid dualizing complexes via differential graded algebras
  (survey)}, Triangulated categories, London Math. Soc. Lecture Note Ser., vol.
  375, Cambridge Univ. Press, Cambridge, 2010, pp.~452--463.

\bibitem{Yekutieli16}
\bysame, \emph{The squaring operation for commutative {DG} rings}, J. Algebra
  \textbf{449} (2016), 50--107.

\bibitem{Yekutieli-Zhang04}
Amnon Yekutieli and James~J. Zhang, \emph{Rigid dualizing complexes on
  schemes}, arXiv math.AG/0405570.

\bibitem{Yekutieli-Zhang97}
\bysame, \emph{Serre duality for noncommutative projective schemes}, Proc.
  Amer. Math. Soc. \textbf{125} (1997), no.~3, 697--707.

\bibitem{Yekutieli-Zhang99}
\bysame, \emph{Rings with {A}uslander dualizing complexes}, J. Algebra
  \textbf{213} (1999), no.~1, 1--51.

\bibitem{Yekutieli-Zhang05}
\bysame, \emph{Dualizing complexes and perverse modules over differential
  algebras}, Compos. Math. \textbf{141} (2005), no.~3, 620--654.

\bibitem{Yekutieli-Zhang06}
\bysame, \emph{Dualizing complexes and perverse sheaves on noncommutative
  ringed schemes}, Selecta Math. (N.S.) \textbf{12} (2006), no.~1, 137--177.

\bibitem{Yekutieli-Zhang08}
\bysame, \emph{Rigid complexes via {DG} algebras}, Trans. Amer. Math. Soc.
  \textbf{360} (2008), no.~6, 3211--3248.

\bibitem{Yekutieli-Zhang09}
\bysame, \emph{Rigid dualizing complexes over commutative rings}, Algebr.
  Represent. Theory \textbf{12} (2009), no.~1, 19--52.

\end{thebibliography}
\end{document}